\documentclass[11pt,reqno,twoside]{amsart}
\usepackage[latin1]{inputenc}
\usepackage{amssymb,amsfonts,amsmath}
\usepackage[english]{babel}
\usepackage{longtable}
\usepackage{graphicx, subfigure}
\usepackage{natbib}
\usepackage[foot]{amsaddr}
\usepackage{array}
\usepackage{color}
\usepackage{enumerate}
\usepackage[breaklinks]{hyperref}

\newtheorem{Theorem}{Theorem}
\newtheorem{Proposition}{Proposition}
\newtheorem{Definition}{Definition}
\newtheorem{Remark}{Remark}
\newtheorem{lemma}{Lemma}
\newenvironment{pf}{\medskip\noindent{\bf Proof.}\enspace}
{\hfill\newline\smallskip}

\begin{document}

\title[A hybrid model of collective motion under alignment and chemotaxis]{A hybrid mathematical model of collective motion under alignment and chemotaxis}

\author[E. Di Costanzo]{Ezio Di Costanzo$^1$}
\address{$^{1,4,5}$Istituto per le Applicazioni del Calcolo ``M. Picone'' -- Consiglio Nazionale delle Ricerche, Italy.}
\author[M. Menci]{Marta Menci$^2$}
\address{$^{2}$Universit\`a Campus Bio-Medico di Roma, Via \`Alvaro del Portillo, 00128, Roma, Italy.}
\author[E. Messina]{Eleonora Messina$^3$}
\address{$^3$Dipartimento di Matematica e Applicazioni ``R. Caccioppoli'', Universit\`a degli Studi di Napoli ``Federico II'', Via Cintia, I-80126 Napoli, Italy.}

\author[R. Natalini]{Roberto Natalini$^4$}

\author[A. Vecchio]{Antonia Vecchio$^5$}

\email{$^1$e.dicostanzo@na.iac.cnr.it}
\email{$^2$m.menci@unicampus.it}
\email{$^3$eleonora.messina@unina.it}
\email{$^4$roberto.natalini@cnr.it}
\email{$^5$antonia.vecchio@cnr.it}



\begin{abstract}
In this paper we propose and study a hybrid \emph{discrete in continuous} mathematical model of collective motion under alignment and chemotaxis effect. Starting from the paper by \citet{dicostanzo}, in which the Cucker-Smale model \citep{cucker} was coupled with other cell mechanisms, to describe the cell migration and self-organization in the zebrafish lateral line primordium, we introduce a simplified model in which the coupling between an alignment and chemotaxis mechanism acts on a system of interacting particles. In particular we rely on a hybrid description in which the agents are discrete entities, while the chemoattractant is considered as a continuous signal. The proposed model is then studied both from an analytical and a numerical point of view. From the analytic point of view we prove, globally in time, existence and uniqueness of the solution. Then, the asymptotic behaviour of a linearised version of the system is investigated. Through a suitable Lyapunov functional we show that for $t\rightarrow +\infty$, the migrating aggregate exponentially converges to a state in which all the particles have a same position with zero velocity. Finally, we present a comparison between the analytical findings and some numerical results, concerning the behaviour of the full nonlinear system.
\end{abstract}
\maketitle
\numberwithin{equation}{section}
\numberwithin{figure}{section}
\section{Introduction}
A collective motion is a form of collective behaviour in which the individual unit's action is strongly dominated by the influence of other units, so that its motion results very different from how it would be if the individual was alone. Collective motion has been extensively studied in recent years in a great variety of systems, from non-living systems, such as nematic fluids, nano swimmers, simple robots, to living systems, such as bacteria colonies, cell aggregates, swarms of insects, flocks of birds, schools of fish, or crowds of humans. Many references can be found, for example, in the review by \citet{vicsek-collective}.
\par Different mechanisms of aggregations in a given system have been proposed. Coordinated motion of cells results in making a biological process more efficient (e.g. in embryogenesis, wound healing, immune response, etc.), while in the case of tumour cell invasion it appears to speed up the progression of the decease \citep{vicsek-collective-cell}. For the animals entities it is observed that a group can more efficiently explore surrounding environments, to enhance foraging capability and detection of predators \citep{pitcher}; being in a group offers protection against attacks by predators \citep{ioannou} and increases the locomotion efficiency \citep{fish}.
\par The research in the field of collective motion modelling is of great interest for the applications in many fields of the real life, from biomedical field, e.g. in relation to the collective motion of cells (see \citealp{szabo,belmonte,arboleda,sepulveda,joie,colin,scianna,dicostanzo}; and the reviews by \citealp{hatzikirou} and \citealp{vicsek-collective-cell}), to socio-economic field and life sciences \citep{pareschi,pareschi2}, to the problems of pedestrian flows \citep{helbing,piccoli,faria,cristiani-animal-pedestrians,bruno,moussaid,tosin-book,cristiani}.
\par Basic mathematical models about collective motion are substantially based on one or more of the following steering behaviours of the units: \emph{alignment}, in order to move in the same direction of average heading direction of neighbouring units (\citealp{vicsek}; \citealp{cucker}, and references therein); \emph{separation}, in order to avoid crowding, and \emph{cohesion}, to remain close to the average position of neighbours \citep{dorsogna,strombom}. Sometimes \emph{three-zone models} combine together alignment, attraction and repulsion effects on three non-overlapping regions, which take into account the sensory capabilities of the individual \citep{aoki,couzin,gregoire1,gregoire2}. For example, the visual field of a bird does not extend behind its body, fishes can accompany visual signals with those coming from their lateral line, while cells can feel stimuli around them. This leads to introduce the concept of the \emph{cone of vision} \citep{huth,couzin,hemelrijk}. In general the cone of interaction is not only animal-dependent, but it can also vary depending on the type of motion, environmental conditions, presence of predators, and aim of the displacement \citep{cristiani-animal-groups}. Moreover, when each unit interacts only with those particles which are closer than a predefined distance, that is a fixed range of influence, we speak about ``metric'' interaction. On the other hand, recent studies of starling flocks, have shown that each bird modifies its position relative to the mates directly surrounding it (typically six or seven individuals), no matter how close or how far away those animals are. In this type of behaviour, referred as ``topological'' interaction \citep{ballerini}, the radius of perception is adjusted by each individual, in such way that the neighbourhood of interaction encompasses a predefined mass of other individuals felt comfortable to interact with it (\citealp{cristiani-animal-pedestrians}; \citealp{cristiani-animal-groups}, and references therein).
\par With regard to the alignment models a widespread model is represented by the Cucker-Smale model \citep{cucker}. The original model was proposed to describe the dynamics in a flock of birds, but its applicability is to general phenomena where autonomous agents reach a consensus, e.g. animal herding, emergence of common languages in primitive societies, etc. \citep{couzin-leader}. The starting point of this model is represented by the seminal paper by \citet{vicsek} and previous analytic studies can be found also in \citet{tsitsiklis} and \citet{jadbabaie}. Anyway, the Cucker-Smale paper is widely known to have established an analytical exact result on the convergence to the same velocity in a group of interacting agents through an alignment effect.
\par The hypothesis of the Cucker-Smale model is that the force acting on every particle (bird) is a weighted average of the differences of its velocity with those of the other particles (birds). In particular, for $N \in \mathbb{N}$ particles, the proposed model, in continuous time $t\in \mathbb{R}_{\geq 0}$, can be written in the form
\begin{align}\label{sys-cuck}
\left\{
	\begin{array}{l}
	\dot{\mathbf{V}}_i=\displaystyle\frac{1}{N}\displaystyle\sum_{j=1}^N\frac{\alpha_1}{\left(\alpha_2+\left\|\mathbf{X}_i-\mathbf{X}_j\right\|^2\right)^{\sigma}}(\mathbf{V}_j-\mathbf{V}_i),\\
	\dot{\mathbf{X}}_i=\mathbf{V}_i,
	\end{array}
	\right.
\end{align}
where $\mathbf{X}_i,\mathbf{V}_i \in\mathbb{R}^n$, for $i=1,\dots,N$, are the position and the velocity of the $i$-th particle, $\alpha_1$, $\alpha_2$, $\sigma$ are positive constants and $\left\|\cdot\right\|$ is the Euclidean norm in $\mathbb{R}^n$. In particular $\sigma$ captures the rate of decay of the influence between agents in the flock as they separate in space, and it is a fundamental parameter in the time-asymptotic behaviour of system \eqref{sys-cuck} (see below). In particular, under suitable initial data and parameters the model implies that the state of the group converges to one in which all particles move with the same velocity, said \emph{flocking state}.
\par The main convergence result proved in \citet{cucker} has been improved by \citet{ha}, using an explicit Lyapunov functional approach. First, as in this last reference, we can give the definition of \emph{time-asymptotic flocking} as follows:
\begin{Definition}\label{def-flocking}
Given system \eqref{sys-cuck}, let $\mathbf{X}_{\text{CM}}=\frac{1}{N}\sum_{i=1}^N \mathbf{X}_i $ and $\mathbf{V}_{\text{CM}}=\frac{1}{N}\sum_{i=1}^N \mathbf{V}_i $ position and velocity of the centre of mass. System \eqref{sys-cuck} has a time-asymptotic flocking if and only if $\left(\mathbf{X}_i,\mathbf{V}_i\right)$, $i=1,\dots,N$, satisfy the two conditions:
\begin{enumerate}[1)]
\item the velocity fluctuations go to zero time-asymptotically (velocity alignment):
\begin{align}\label{lim-varianza-x}
	\lim_{t \rightarrow +\infty}\sum_{i=1}^N \left\|\mathbf{V}_i(t)-\mathbf{V}_\text{CM}(t)\right\|^2=0;
	\end{align}
	\item the position fluctuations are uniformly bounded in time $t$ (forming a group):
	\begin{align}\label{lim-varianza-v}
	\sup_{0\leq t\leq +\infty}\sum_{i=1}^N \left\|\mathbf{X}_i(t)-\mathbf{X}_\text{CM}(t)\right\|^2<+\infty.
	\end{align}
\end{enumerate}
\end{Definition}
Notice that the square root of the quantities under the limit and supremum operations in \eqref{lim-varianza-x}--\eqref{lim-varianza-v} is proportional to the standard deviations of $\mathbf{V}_i(t)$ and $\mathbf{X}_i(t)$ around the centre of mass system.
\citet{cucker} and \citet{ha} proved that for $\sigma \in[0,1/2]$ occurs a global \emph{unconditional flocking} of system \eqref{sys-cuck}, as stated in Definition \ref{def-flocking}, regardless of initial configurations, while for $\sigma \in(1/2,+\infty)$ there is \emph{conditional flocking}, that is only some parameters and initial data lead to a flocking state, but in general the dispersion of the flock may occur.
\par Cucker and Smale-like models have been widely employed in collective dynamics and several applications can be found, from the biological field \citep{szabo,belmonte,arboleda,sepulveda,dicostanzo,dicostanzoITM,doi:10.1093/imammb/dqw022}, to the collective motion of different interacting groups \citep{albi}. Furthermore, many extensions have been proposed. For example, \citet{cucker-dong,cucker-dong1} added a repelling force between particles in equation \eqref{sys-cuck}$_1$, proving an analogous convergence theorem ensuring, on certain conditions, flocking behaviour, in addition without collision between particles. \citet{ha2} applied system \eqref{sys-cuck} to describe the motion of phototactic bacteria (i.e., bacteria that move towards light), adding a force that move excited bacteria towards the light source in equation \eqref{sys-cuck}$_1$, and an additional equation for the excitation level of each particle. Also in this case, under particular conditions, it is proved that the asymptotic velocity of the particles tends to an identical terminal velocity. Other models have introduced noise \citep{cucker-mordecki}, stochastic equations \citep{ha-lee}, leader individuals with a preferred heading direction \citep{cucker-huepe}, kinetic equations \citep{carrillo}.
\par In this paper, starting from the model proposed by \citet{dicostanzo} in relation to the morphogenesis in the zebrafish lateral line primordium, in which the authors coupled the Cucker-Smale model with other cell mechanisms (chemotaxis, attraction-repulsion, damping effects), we study, both from an analytical and a numerical point of view, a simplified model in which the coupling between an alignment and chemotaxis effect acts in a system of interacting particles. Our description is hybrid: discrete for the interacting agents, and continuous for the chemotactic signal. From the analytic point of view, in the two-dimensional case with $N$ particles, we prove, globally in time, existence and uniqueness of solutions to the model. Then, the time-asymptotic behaviour of the model is investigated. In particular, through a suitable Lyapunov functional, it is shown that position and velocity of all the particles go exponentially to those of their centre of mass. Moreover, the velocity of the centre of mass tends time-asymptotically to zero. From a numerical point of view, on a bounded domain with periodic boundary conditions, some 2D dynamical tests are proposed to simulate the behaviour of the model, and to compare the numerical to the analytical results.
\par From our findings, we observe that while in the Cucker-Smale model two conditions may occur, conditional and unconditional flocking, in our model the flocking behaviour, as given in Definition \ref{def-flocking}, is always ensured. Moreover, we have the stronger conditions that all the particles converge in a single position, i.e. their centre of mass, and that the velocity of the centre of mass tends to zero. The result is somewhat as expected, but it can be considered a first rigorous study of these hybrid model.
\par Among the various numerical simulations, we discuss the competition between alignment and chemotactic effects varying the parameters of the system. We find a decrease in the rate of convergence of the particles when the strength of the alignment term increases with respect to the chemotaxis. On the other hand, an increase in the rate of convergence can be found with the same parameters and with a greater number of interacting particles. Inspired then from the aforesaid paper on the zebrafish lateral line, we consider also the case of two kinds of cells: leader cells corresponding to the sources of the chemotactic signal, and follower cells that do not produce any chemical signal, both subject to the alignment and chemotaxis force. In this case we observe the convergence of the group toward the sources of the chemoattractant. Finally, we consider the collective motion under chemotaxis, neglecting the alignment effect. Our numerical tests, in this regard, show the absence of convergence, and an oscillating motion of the particles around their centre of mass. This suggests that, in such model, the only chemotactic effect is unable to reproduce biological phenomena involving stationary aggregates.
\par The paper is organized as follows. In Section \ref{sec:basic-flock-chemo} we design our hybrid mathematical model of collective motion. Section \ref{sec:local} deals with the analytical study of the local existence and uniqueness of the solution. Then, in Section \ref{sec:global} this result is extended and the global existence and uniqueness is proved. Section \ref{sec-stab} is devoted to investigate the asymptotic behaviour of the model around the equilibrium configurations. In Section \ref{sec:numerical-chemo} we show and discuss some numerical results. Finally, Section \ref{sec:conclusion-flock} includes the conclusions and possible future perspectives.
\section{The basic mathematical model}\label{sec:basic-flock-chemo}
Starting from model \eqref{sys-cuck}, we assume that the force acting on each particle is given by an alignment term, proportional to the differences of velocity with the other particles and weighed on the distances, and by a chemotactic attraction towards higher concentration of a chemical signal $f(\mathbf{x},t)$, supposed produced by the particles themselves. Typically, this last force is proportional to the gradient of the concentration $\nabla f$ (see \citealp{eisenbach} for biological backgrounds, and \citealp{murrayII,perthame} for some mathematical references). In our hybrid description, while particles are considered discrete entities, endowed of a radius $R$ describing their circular shape, the signal $f$ is supposed to be continuous and its rate of change in time is equal to a diffusion term, a source term depending on the position of each particle, and a degradation term.

\par To summarize our hypotheses, we write the following system:
\begin{align}\label{sys-cuck-chemo}
\left\{
	\begin{array}{l}
	\dot{\mathbf{V}}_i=\displaystyle\frac{\beta}{N}\displaystyle\sum_{j=1}^N\frac{1}{\left(1+\frac{\left\|\mathbf{X}_i-\mathbf{X}_j\right\|^2}{R^2}\right)^{\sigma}}(\mathbf{V}_j-\mathbf{V}_i)+\gamma \nabla f (\mathbf{X}_i),\\
	\dot{\mathbf{X}}_i=\mathbf{V}_i,\\\\
	\partial_t f=D\Delta f+ \displaystyle g_{\mathbf{X}}-\eta f,
	\end{array}
	\right.
\end{align}
Initial data are given by initial position and velocity for each particle:
\begin{align*}
	\mathbf{X}(0)&=\mathbf{X}_0,\quad 	 \mathbf{V}(0)=\mathbf{V}_0,
\end{align*}
with $\mathbf{X}=\left(\mathbf{X}_1,\dots,\mathbf{X}_N\right)$, $\mathbf{V}=\left(\mathbf{V}_1,\dots,\mathbf{V}_N\right)$, and by the initial concentration of signal, that we assume
\begin{align}\label{iniz-f}
	f(\mathbf{x},0):=f_0=0.
\end{align}
We note that equation \eqref{sys-cuck-chemo}$_3$ can be analytically solved making the classical exponential transformation:
\begin{align*}
	f=e^{-\eta t} u,
\end{align*}
with $u(\mathbf{x},t)$ solution of
\begin{align*}
	\partial_t u=D\Delta u+ e^{\eta t} g_{\mathbf{X}},
\end{align*}
and $u(\mathbf{x},0)=f_0$.
Now, if
\begin{align*}
	\Gamma(\mathbf{x},t):=\frac{1}{(4\pi D t)^{n/2}}e^{-\frac{||\mathbf{x}||^2}{4Dt}}
\end{align*}
is the fundamental solution of the heat equation on $\mathbb{R}^n$, the unique solution of \eqref{sys-cuck-chemo}$_3$ can be written
\begin{align}
	f(\mathbf{x},t)&=(\Gamma(\mathbf{x},t) \ast f_0)e^{-\eta t}+e^{-\eta t}\int_0^t\Gamma(\mathbf{x},t-\tau)\ast \left(e^{\eta \tau}g_{\mathbf{X}}\left(  \mathbf{x},\tau \right)\right)\,d \tau\nonumber\\
	&= \int_0^t e^{-\eta (t-\tau)} \int_{\mathbb{R}^n}\Gamma(\mathbf{x}-\bar{\mathbf{x}},t-\tau)g_{\mathbf{X}}\left(  \overline{\mathbf{x}},\tau \right)\,d \bar{\mathbf{x}} d \tau, \label{soluz-fond}
\end{align}
where $\ast$ is the convolution operation in the variable $\mathbf{x}$, and $f_0=0$ for our initial condition \eqref{iniz-f}.
\par In the following, for analytical and numerical simplicity, we will discuss the case of $N$ particles in $\mathbb{R}^2$.
\par For the source term $g_{\mathbf{X}}$ in the chemical signal equation,  we get
\begin{align}\label{scelta_g}
	g_{\mathbf{X}}(\mathbf{x},t):=\xi \sum_{j=1}^{N} \varphi\left(\left| x_1-X_{j1} \right|\right)\varphi\left(\left| x_2-X_{j2} \right|\right),
\end{align}
where $\xi>0$, $\mathbf{x}=\left(x_{1},x_{2}\right)$, $\mathbf{X}_j=\left(X_{j1},X_{j2}\right)$ is the position of $j-$th cell, and $\varphi: \mathbb{R}\rightarrow \mathbb{R}$ is a nonnegative bounded function, $\varphi \in C^{1}_{c}\left(\left[-R,R\right]\right)$. This function is intended to take into account for the production of chemical signal due to each cell. 

 First, if $\mathbf{x}=(x_1,x_2)$ and $\bar{\mathbf{x}}=(\bar{x}_1,\bar{x}_2)$, from \eqref{soluz-fond} and \eqref{scelta_g} we can thus write
\begin{align}\label{calore2D}
	f(x_1,x_2,t)=& \xi\sum_{j=1}^N\int_0^t \frac{e^{-\eta (t-\tau)}}{4\pi(t-\tau)D}
	\int_{\mathbb{R}} \phi\left( x_1-\overline{x}_1\right) \varphi\left( \overline{x}_1-X_{j1}(\tau)\right)\,d\bar{x}_1 \nonumber \\
	 &\cdot \int_{\mathbb{R}} \phi\left( x_2-\overline{x}_2\right) \varphi\left( \overline{x}_2-X_{j2}(\tau)\right)\,d\bar{x}_2\,d\tau,
\end{align}
where, for the sake of notational simplicity, we denote with $\phi:\mathbb{R} \times (0,T] \longrightarrow \mathbb{R}$ the function given by
$\phi\left(x,t\right)=e^{-\frac{x^2}{4Dt}}$.\\
Equation \eqref{calore2D} can be rewritten in a more compact way, introducing \mbox{$W: E \longrightarrow C\left(\mathbb{R} \times \left(0,T \right]\right)$}, with \mbox{$E:= \left\{ \psi \in C\left(\mathbb{R}\right) : \left| \psi(x) \right| \le Ce^{hx^2} \ \ \mbox{with} \ \  C>0, h<\frac{1}{4DT}\right\}$}, defined by

\begin{equation}\label{defW}
W\left( \psi \right) \left(x,t \right):= \int_{\mathbb{R}} \frac{1}{\sqrt{4\pi D t} }\phi\left(x-y,t\right)\psi\left(y\right) dy.
\end{equation}
We observe that, \eqref{defW} is the unique classical solution to the parabolic Cauchy problem, given by
\begin{align*}
\left\{
\begin{array}{lll}
\partial_t W(\psi)(x,t)=D\partial_{xx}^{2} W(\psi)(x,t), \quad (x,t) \in \mathbb{R}\times (0,T], \\\\
W(\psi)(x,0)=\psi(x).
\end{array}
\right.
\end{align*}
Hence, \eqref{calore2D} can be rewritten as
 \begin{align}
	f(x_1,x_2,t)= \xi\sum_{j=1}^N\int_0^t e^{-\eta (t-\tau)}
	W(\varphi) \left( x_1-X_{j1}(\tau), t-\tau \right) W(\varphi) \left( x_2-X_{j2}(\tau), t-\tau \right)\,d\tau.
\end{align}
Before computing the chemotactic gradient $\nabla f=(\partial_{x_1} f,\partial_{x_2} f)$, we observe that
for any $\psi \in E$, $(x,t) \in \mathbb{R}\times (0,T]$, it follows
\begin{align*}
\partial_x W(\psi)(x,t)= \int_{\mathbb{R}} \frac{1}{\sqrt{4\pi D t}} \left(-\frac{x-y}{2Dt} \phi\left(x-y,t\right) \right) \psi\left(y\right) dy =  \int_{\mathbb{R}} \frac{\phi\left(x-y,t\right) \psi^{\prime} \left(y\right)}{\sqrt{4\pi D t}} dy=W(\psi^{\prime})(x,t).
\end{align*}
Thus 
\begin{align}\label{gradientef1}
\partial_{x_1}  f(x_1,x_2,t)= \xi\sum_{j=1}^N\int_0^t e^{-\eta (t-\tau)}
	W(\varphi^{\prime}) \left( x_1-X_{j1}(\tau), t-\tau \right) W(\varphi) \left( x_2-X_{j2}(\tau), t-\tau \right)\,d\tau, \nonumber\\
\partial_{x_2}  f(x_1,x_2,t)= \xi\sum_{j=1}^N\int_0^t e^{-\eta (t-\tau)}
	W(\varphi) \left( x_1-X_{j1}(\tau), t-\tau \right) W(\varphi^{\prime}) \left( x_2-X_{j2}(\tau), t-\tau \right)\,d\tau.
\end{align}
\par Finally, substituting \eqref{gradientef1} into \eqref{sys-cuck-chemo}$_1$ we can summarize, for $i=1,\dots,N$, the following system:
\begin{align}\label{Nparticles2D_new}
\left\{
\begin{array}{ll}
	\dot{V}_{i1}&=\displaystyle\frac{\beta}{N}\sum_{j=1}^N\frac{1}{\left(1+\frac{||\mathbf{X}_i-\mathbf{X}_j||^2}{R^2}\right)^{\sigma}}(V_{j1}-V_{i1})\\
	&+  \displaystyle \sum_{j=1}^{N}  \int_{0}^{t} K(t-\tau)  W(\varphi^{\prime}) \left( X_{i1}(t)-X_{j1}(\tau), t-\tau \right) W(\varphi) \left( X_{i2}(t)-X_{j2}(\tau), t-\tau \right)\,d\tau,\\	
	\dot{V}_{i2}&=\displaystyle\frac{\beta}{N}\sum_{j=1}^N\frac{1}{\left(1+\frac{||\mathbf{X}_i-\mathbf{X}_j||^2}{R^2}\right)^{\sigma}}(V_{j2}-V_{i2})\\
	&+ \displaystyle \sum_{j=1}^{N}  \int_{0}^{t} K(t-\tau) W(\varphi) \left(  X_{i1}(t)-X_{j1}(\tau), t-\tau \right) W(\varphi^{'}) \left( X_{i2}(t)-X_{j2}(\tau), t-\tau \right)\,d\tau,\\
	\dot{X}_{i1}&=V_{i1},\\
	\dot{X}_{i2}&=V_{i2},
	\end{array}
	\right.
\end{align}
with
	\begin{align}\label{K-stab}
K(t-\tau):=\gamma\xi e^{-\eta (t-\tau)}.
	\end{align}

\section{Local existence and uniqueness of the solution}\label{sec:local}
In this section, using a fixed point argumentation, we prove for \eqref{Nparticles2D_new} the local existence and uniqueness of solutions. In the next section, then, this result will be extended to a global result in time.
\par First, let $\mathbf{y}=(\mathbf{V}_1,\dots,\mathbf{V}_N,\mathbf{X}_1,\dots,\mathbf{X}_N)$ the solution vector, and let
\begin{align*}
	\mathbf{q}:=(\mathbf{q}_{1},\mathbf{q}_{2}),
\end{align*}
with
\begin{align}\label{qi}
	q_{1,i1}&:=\frac{\beta}{N}\sum_{j=1}^N\frac{1}{\left(1+\frac{\left\|\mathbf{X}_j-\mathbf{X}_i\right\|^2}{R^2}\right)^{\sigma}}(V_{j1}-V_{i1}),\\
	q_{1,i2}&:=\frac{\beta}{N}\sum_{j=1}^N\frac{1}{\left(1+\frac{\left\|\mathbf{X}_j-\mathbf{X}_i\right\|^2}{R^2}\right)^{\sigma}}(V_{j2}-V_{i2}),\\
	q_{2,i1}&:=V_{i1},\\
	q_{2,i2}&:=V_{i2}.\label{qi_end}
\end{align}
and then
\begin{align*}
\mathbf{p}:=(\mathbf{p}_{1},\mathbf{p}_{2}),
\end{align*}
with
	\begin{align}
	p_{1,i1}&:=\sum_{j=1}^N W(\varphi^{'}) \left( X_{i1}(t)-X_{j1}(\tau), t-\tau \right) W(\varphi) \left( X_{i2}(t)-X_{j2}(\tau), t-\tau \right) \label{p1-i1},\\
	p_{1,i2}&:=\sum_{j=1}^N W(\varphi) \left( X_{i1}(t)-X_{j1}(\tau), t-\tau \right) W(\varphi^{'}) \left( X_{i2}(t)-X_{j2}(\tau), t-\tau \right)\label{p1-i2},\\
	p_{2,i1}&=p_{2,i2}=0.\nonumber
\end{align}
Finally, 
system \eqref{Nparticles2D_new} can be written as
\begin{align*}
	\dot{\mathbf{y}}=\mathbf{q}(\mathbf{y})+\int_0^t K(t-\tau) \mathbf{p}(t-\tau,\mathbf{y}(t),\mathbf{y}(\tau))\,d\tau.
\end{align*}
with $K(t-\tau)$ given in \eqref{K-stab}. Integrating from 0 to $t$ we have
\begin{align*}
	\mathbf{y}=\mathbf{y}_0+\int_0^t \mathbf{q}(\mathbf{y}(\tau))\,d\tau+\int_0^t\int_0^s K(s-\tau)\mathbf{p}(s-\tau,\mathbf{y}(s),\mathbf{y}(\tau))\,d\tau\,ds,
\end{align*}
with $\mathbf{y}_0=\mathbf{y}(0)$. Then, interchanging the order of integration in the second integral, we have
\begin{align}\label{eq-integrale1}
	\mathbf{y}=\mathbf{y}_0+\int_0^t\left[ \mathbf{q}(\mathbf{y}(\tau))+\int_\tau^t K(s-\tau) \mathbf{p}(s-\tau,\mathbf{y}(s),\mathbf{y}(\tau))\,d s\right]\,d\tau,
\end{align}
or
\begin{align}\label{eq-integrale2}
	\mathbf{y}=\mathbf{y}_0+\int_0^t\left[ \mathbf{q}(\mathbf{y}(\tau))+\mathbf{h}(t,\tau,\mathbf{y}(\tau))\right]\,d\tau,
\end{align}
with
\begin{align*}
	\mathbf{h}(t,\tau,\mathbf{y}(\tau)):=\int_\tau^t K(s-\tau) \mathbf{p}(s-\tau,\mathbf{y}(s),\mathbf{y}(\tau))\,ds.
\end{align*}
For a discussion of such type of equations see, for example, \citet{burton}, \citet{khalil}, and also \citet{wazwaz}, \citet{lakshmikantham}.
\par Now, let $a,b>0$. We consider the set
\begin{align*}
	S=\left\{(t,\tau,s,\mathbf{y}):0\leq \tau \leq s< t\leq a,\,\left\|\mathbf{y}(t)-\mathbf{y}_0\right\|\leq b\right\}.
\end{align*}
Since $\mathbf{q}(\mathbf{y})$ is continuous on $S$, we can define
\begin{align}\label{M1}
	M_1=\max_S \left\| \mathbf{q}(\mathbf{y})\right\|.
\end{align}
Then, to prove that $\mathbf{h}(t,\tau,\mathbf{y}(\tau))$ is continuous in $S$, firstly we prove that $K(s-\tau)\mathbf{p}(s-\tau,\mathbf{y}(s),\mathbf{y}(\tau))$ in $L^1(\tau,t)$ with respect to variable $s$. It is enough to demonstrate the integrability around $s=\tau$. 
Starting from \eqref{p1-i1} and \eqref{defW}, we consider $p_{1,i1}(s-\tau,\mathbf{y}(s),\mathbf{y}(\tau))$, and the change of variables
\begin{align}\label{cambio1}
	\frac{X_{i1}-X_{j1}-\tilde{y}_1}{\sqrt{4(s-\tau)D}}=z_1, \ \ \ \ \ \ \ \ \ \ \ \ \ \frac{X_{i2}-X_{j2}-\tilde{y}_2}{\sqrt{4(s-\tau)D}}=z_2.
\end{align}

It follows
\begin{align*}
	\left|p_{1,i1}(s-\tau,\mathbf{y}(s),\mathbf{y}(\tau))\right|&\leq  \frac{\sqrt{4D(s-\tau)}\sqrt{4D(s-\tau)}}{4D(s-\tau)\pi}\sum_{j=1}^N \int_{\mathbb{R}}e^{-z_1^2}dz_1  \int_{\mathbb{R}}e^{-z_2^2}dz_2 \left| \left| \varphi\right| \right|  \left| \left| \varphi^{\prime}\right| \right|\\
	&= N \left| \left| \varphi\right| \right|  \left| \left| \varphi^\prime\right| \right|.
	\end{align*}
The same holds for $p_{1,i2}$ in \eqref{p1-i2}, so we can write

\begin{align}\label{cont-Cp}
 \left\|\mathbf{p}(s-\tau,\mathbf{y}(s),\mathbf{y}(\tau))\right\| \leq \sqrt{2N^3} \left| \left| \varphi \right| \right|  \left| \left| \varphi^\prime\right| \right|
\end{align}
and

\begin{align}\label{cont-Cp2}
\left|K(s-\tau)\right|  \left\|\mathbf{p}(s-\tau,\mathbf{y}(s),\mathbf{y}(\tau)) \right\|\leq \gamma \xi e^{-\eta(s-\tau)}	\sqrt{2N^3} \left| \left| \varphi\right| \right|  \left| \left| \varphi^\prime\right| \right| .
\end{align}
Now $K(s-\tau)\mathbf{p}(s-\tau,\mathbf{y}(s),\mathbf{y}(\tau))$ is continuous in $\mathbf{y}$ and, from \eqref{cont-Cp2}, it is $L^1(\tau,t)$ with respect to the variable $s$, so $\mathbf{h}(t,\tau,\mathbf{y}(\tau))$ is continuous in $S$, and we can define
\begin{align}\label{M2}
	M_2=\max_S \left\|\mathbf{h}(t,\tau,\mathbf{y}(\tau)) \right\|.
\end{align}

To prove local existence and uniqueness, we want to obtain a Lipschitz condition in $S$ for the functions $\mathbf{q}$ and $\mathbf{p}$ with respect to the variable $\mathbf{y}$.
First, because of $\mathbf{q}$ is $C^1$ on $S$, the Jacobian matrix $\left[\partial \mathbf{q}/\partial \mathbf{y}\right]$ is bounded on $S$ uniformly in $\tau$,
so $\mathbf{q}$ satisfies the Lipschitz condition
\begin{align}\label{cond-lip-q}
	\left\|\mathbf{q}(\mathbf{y}_1)-\mathbf{q}(\mathbf{y}_2)\right\|\leq L_1 \left\| \mathbf{y}_1-\mathbf{y}_2\right\|,
\end{align}
with $L_1$ positive constant and $(t,\tau,s,\mathbf{y}_1)$, $(t,\tau,s,\mathbf{y}_2)\in S$.\\
To obtain a Lipschitz condition in S for $\textbf{p}$ with respect to $\textbf{y}$, we preliminary observe that, for any $\psi \in E$, $\left( x,t\right), \left( x^\prime,t\right) \in \mathbb{R} \times \left(0,T\right]$, by definition (\ref{defW}) immediately follows

\begin{equation}
\left| W\left( \psi \right) (x,t) \right| \leq \frac{||\psi||}{\sqrt{\pi}} \int_{\mathbb{R}} e^{-z^2}dz= ||\psi||.
\end{equation}

Moreover, using the mean value theorem, and the estimate of the derivative of function $\Gamma$ \citep{friedman2008partial}, we obtain
\begin{align}
\left| W\left( \psi \right) (x,t)-W\left( \psi \right) (x^\prime,t) \right|& \leq 
 \int_{\mathbb{R}} \left| \Gamma\left(x-y,t\right)- \Gamma \left(x^\prime-y,t\right)  \right| \left| \psi \left( y\right) \right| dy  \\ \nonumber
& \leq  \left| \left| \psi \right| \right| \frac{\left| x-x^\prime\right|}{t} \int_{\mathbb{R}}e^{-|u|^2}\sqrt{4Dt} \ du \\ \nonumber
& \leq  2 \sqrt{D\pi} \left| \left| \psi \right| \right| \frac{\left| x-x^\prime\right|}{\sqrt{t}} . 
\end{align}
Let now focus on the component $p_{1,i1}$ of $\textbf{p}$. We observe that, for any $\left(t,\tau,s,\mathbf{y_1} \right)$, $\left(t,\tau,s,\mathbf{y_2} \right) \in S$, denoting with $X^{(1)}$, $X^{(2)}$ the variables belonging respectively to $\mathbf{y_1}$ and $\mathbf{y_2}$, we obtain  \\
\begin{align}
&\left|p_{1,i1}(s-\tau,\mathbf{y}_1(s),\mathbf{y}_1(\tau))-p_{1,i1}(s-\tau,\mathbf{y}_2(s),\mathbf{y}_2(\tau))\right| \nonumber \\ \nonumber
&\leq\sum_{j=1}^{N} \left| W(\varphi^\prime) \left(X^{(1)}_{i1}\left(s\right)-X^{(1)}_{j1}(\tau), s-\tau \right) W(\varphi) \left( X^{(1)}_{i2}\left(s\right)-X^{(1)}_{j2}(\tau), s-\tau \right) \right.\\ \nonumber
&\left. \ \ \    - W(\varphi^\prime) \left(X^{(2)}_{i1}\left(s\right)-X^{(2)}_{j1}(\tau), s-\tau \right) W(\varphi) \left( X^{(2)}_{i2}\left(s\right)-X^{(2)}_{j2}(\tau), s-\tau \right) \right| \\  \nonumber
&\leq \sum_{j=1}^{N} \left(\left|  W(\varphi^\prime) \left(X^{(1)}_{i1}\left(s\right)-X^{(1)}_{j1}(\tau), s-\tau \right) \right|
\left| W(\varphi) \left( X^{(1)}_{i2}\left(s\right)-X^{(1)}_{j2}(\tau), s-\tau \right) \right.\right. \\ \nonumber
&\left. \ \ \    -W(\varphi) \left( X^{(2)}_{i2}\left(s\right)-X^{(2)}_{j2}(\tau), s-\tau \right) \right| +
\left|  W(\varphi) \left(X^{(2)}_{i2}\left(s\right)-X^{(2)}_{j2}(\tau), s-\tau \right) \right| \\\ \nonumber
& \ \ \    \cdot\left.\left|-W(\varphi^\prime) \left( X^{(2)}_{i1}\left(s\right)-X^{(2)}_{j1}(\tau), s-\tau \right)+ W(\varphi^\prime) \left( X^{(1)}_{i1}\left(s\right)-X^{(1)}_{j1}(\tau), s-\tau \right) \right|\right) \\\ \nonumber
&  \leq\sum_{j=1}^{N} \frac{ 2\sqrt{D\pi} \left|  \left| \varphi \right| \right|  \left|  \left| \varphi^{\prime} \right| \right|}{\sqrt{s-\tau}} \left( \left| X^{(1)}_{i2}\left(s\right)-X^{(2)}_{i2}\left(s\right)\right| +\left| X^{(1)}_{j2}\left(\tau\right)-X^{(2)}_{j2}\left(\tau\right)\right|  \right.  \\ 
& \ \ \   \left. +\left| X^{(2)}_{i1}\left(s\right)-X^{(1)}_{i1}\left(s\right)\right| + 
\left| X^{(2)}_{j1}\left(\tau\right)-X^{(1)}_{j1}\left(\tau\right)\right|  \right).  
\end{align}
The same can be done for $p_{1,i2}$, so  $\mathbf{p}$ satisfy the following condition in $\mathbf{y}$ on $S$:
\begin{equation}
	\left\|\mathbf{p}(s-\tau,\mathbf{y}_1(s),\mathbf{y}_1(\tau))-\mathbf{p}(s-\tau,\mathbf{y}_2(s),\mathbf{y}_2(\tau))\right\|
	\leq \frac{L_2}{ \sqrt{s-\tau}}\left(\left\| \mathbf{y}_1(s)-\mathbf{y}_2(s)\right\|
	+\left\| \mathbf{y}_1(\tau)-\mathbf{y}_2(\tau)\right\| \right), \label{cond-lip-p}
	\end{equation}
with $L_2$ a suitable positive constant that incorporates previous constants, and $(\tau,s,\mathbf{y}_1)$, $(\tau,s,\mathbf{y}_2)\in S$.

\par Now, we fix
\begin{align}\label{T-esit-unic}
	T=\min\left[a,\frac{b}{M_1+M_2},\frac{1}{L_1+2 L_2 M}\right],
\end{align}
with $M_1$, $M_2$ given by \eqref{M1}, \eqref{M2}, and
\begin{align}\label{int-infinito-C}
	M:=\int_0^{+\infty}\frac{\left|K(z)\right|}{\sqrt{z}}\,dz.
\end{align}
%
%
%
Then we prove the following
%
%
\begin{Theorem}\label{Theo-ex-un}
Equation \eqref{eq-integrale2} has a unique solution on $\left[0,T\right]$, where $T$ is defined in \eqref{T-esit-unic}.
\end{Theorem}
\begin{pf} 
We consider the functional space
\begin{align*}
	\mathcal{B}=\left\{\mathbf{y}\in C^0\left([0,T]\right):||\mathbf{y}-\mathbf{y}_0||_{C^0}\leq b\right\},
\end{align*}
where
\begin{align*}
	||\mathbf{y}-\mathbf{z}||_{C^0}:=\sup_{0\leq t \leq T}\left\|\mathbf{y}(t)-\mathbf{z}(t)\right\|,
\end{align*}
and we define the functional $\mathbf{A}:\mathcal{B}\rightarrow \mathcal{B}$ as
\begin{align*}
	\mathbf{A}(\mathbf{y})(t):=\mathbf{y}_0+\int_0^t\left[ \mathbf{q}(\mathbf{y}(\tau))+\mathbf{h}(t,\tau,\mathbf{y}(\tau))\right]\,d\tau,
\end{align*}
To see that $\mathbf{A}:\mathcal{B}\rightarrow \mathcal{B}$ notice that $\mathbf{y}$ continuous implies $\mathbf{A}(\mathbf{y})$ continuous, because $\mathbf{q}$ and $\mathbf{h}$ are continuous, and that
\begin{align*}
	\left\|\mathbf{A}(\mathbf{y})-\mathbf{y}_0\right\|_{C^0}&=\sup_{0\leq t \leq T} \left\|\mathbf{A}(\mathbf{y})(t)-\mathbf{y}_0 \right\|\\
	&\leq \sup_{0\leq t \leq T}\int_0^t\left( \left\|\mathbf{q}(\mathbf{y}(\tau))\right\|+\left\|\mathbf{h}(t,\tau,\mathbf{y}(\tau))\right\|\right)\,d\tau\\
	&\leq (M_1 +M_2)T\\
	&\leq b,
\end{align*}
where we have used \eqref{M1}, \eqref{M2} and, in the last inequality, \eqref{T-esit-unic}. To see that $\mathbf{A}$ is a contraction mapping, notice that if $\mathbf{y}_1$ and $\mathbf{y}_2\in \mathcal{B}$ then
\begin{align*}
	\left\|\mathbf{A}(\mathbf{y}_1)-\mathbf{A}(\mathbf{y}_2)\right\|_{C^0}&=\sup_{0\leq t \leq T} \left\|\mathbf{A}(\mathbf{y}_1)(t)-\mathbf{A}(\mathbf{y}_2)(t) \right\|\\
		&\leq \sup_{0\leq t \leq T}\int_0^t( \left\|\mathbf{q}(\mathbf{y}_1(\tau))-\mathbf{q}(\mathbf{y}_2(\tau))\right\|\\
		&+\int_\tau^t\left|K(s-\tau)\right|\left\|\mathbf{p}(s-\tau,\mathbf{y}_1(s),\mathbf{y}_1(\tau))-\mathbf{p}(s-\tau,\mathbf{y}_2(s),\mathbf{y}_2(\tau))\right\|)\,ds\,d\tau\\
		&\leq_{\scriptscriptstyle \eqref{cond-lip-q},\eqref{cond-lip-p}} \sup_{0\leq t \leq T}\int_0^t[ L_1\left\|\mathbf{y}_1(\tau)-\mathbf{y}_2(\tau)\right\|+L_2\int_\tau^t  \frac{1}{\sqrt{s-\tau}}\left|K(s-\tau)\right|\\
		&\cdot \left(\left\|\mathbf{y}_1(s)-\mathbf{y}_2(s)\right\|+\left\|\mathbf{y}_1(\tau)-\mathbf{y}_2(\tau)\right\|\right)]\,ds\,d\tau\\
		&\leq L_1T\left\|\mathbf{y}_1-\mathbf{y}_2\right\|_{C^0} + 2L_2\left\|\mathbf{y}_1-\mathbf{y}_2\right\|_{C^0}\\
		&\cdot\sup_{0\leq t \leq T}\int_0^t\int_\tau^t \left|K(s-\tau)\right|\frac{1}{\sqrt{s-\tau}}\,ds\,d\tau\\
		&= \left(L_1T + 2L_2\sup_{0\leq t \leq T}\int_0^t\int_0^{t-\tau} \left|K(z)\right|\frac{1}{\sqrt{z}}\,dz\,d\tau\right)\left\|\mathbf{y}_1-\mathbf{y}_2\right\|_{C^0}\\
		&\leq\left( L_1T+ 2L_2\int_0^T\int_0^{+\infty} \left|K(z)\right|\frac{1}{\sqrt{z}}\,dz\,d\tau\right)\left\|\mathbf{y}_1-\mathbf{y}_2\right\|_{C^0}\\
		&=( L_1+2L_2 M)T\left\|\mathbf{y}_1-\mathbf{y}_2\right\|_{C^0}.
\end{align*}
From \eqref{T-esit-unic} the constant $( L_1+2L_2 M)T\in (0,1)$. The Banach-Caccioppoli fixed-point theorem completes the proof.
\end{pf}
\section{Global existence of the solution}\label{sec:global}
To obtain global existence for \eqref{eq-integrale1} we will use a principle of continuation of solutions. We will prove that bounded solutions can be continued to $t=+\infty$.
The following general result, adapted to equation \eqref{eq-integrale1}, provides a condition for the continuation of solutions.
\begin{Proposition}\label{princ-prolung}
Let $\mathbf{y}(t)$ be a solution of \eqref{eq-integrale1} on a interval $[0,T)$, if there is a constant $P$ with $\left\|\mathbf{y}-\mathbf{y}_0\right\|\leq P$ on $[0,T)$, then there is a $\bar{T}>T$ such that $\mathbf{y}(t)$ can be continued to $[0,\bar{T}]$.
\end{Proposition}
\begin{pf} We show that $\lim_{t \rightarrow T^{-}} \mathbf{y}(t)$ exists, so we can applied Theorem \ref{Theo-ex-un} starting at $t=T$, and this completes the proof.
\par Let $t_n$ be a monotonic increasing sequence with limit $T$, and let
\begin{align*}
	\bar{U}=\left\{(t,\tau,s,\mathbf{y}):0\leq \tau\leq s\leq t\leq T,\; \left\|\mathbf{y}-\mathbf{y}_0\right\|\leq P\right\}.
\end{align*}
We prove that $\left\{\mathbf{y}(t_n)\right\} $ is a Cauchy sequence. If $t_m>t_n$, from \eqref{eq-integrale1} we have
\begin{align*}
	\left\|\mathbf{y}(t_m)-\mathbf{y}(t_n)\right\|&=\left\|\int_0^{t_m} \left[\mathbf{q}(\mathbf{y}(\tau))+\int_\tau^{t_m} K(s-\tau)\mathbf{p}(s-\tau,\mathbf{y}(s),\mathbf{y}(\tau))\,ds\right]\,d\tau \right.\\
&-\left. \int_0^{t_n} \left[\mathbf{q}(\mathbf{y}(\tau))+\int_\tau^{t_n} K(s-\tau)\mathbf{p}(s-\tau,\mathbf{y}(s),\mathbf{y}(\tau))\,ds\right]\,d\tau\right\|\\
&\leq \int_0^{t_n} \left\| \int_\tau^{t_m} K(s-\tau)\mathbf{p}(s-\tau,\mathbf{y}(s),\mathbf{y}(\tau))\,ds \right.\\
&-\left.\int_\tau^{t_n} K(s-\tau)\mathbf{p}(s-\tau,\mathbf{y}(s),\mathbf{y}(\tau))\,ds \right\|\,d\tau\\
&+\left\|\int_{t_n}^{t_m} \left[\mathbf{q}(\mathbf{y}(\tau))+\int_\tau^{t_m} K(s-\tau)\mathbf{p}(s-\tau,\mathbf{y}(s),\mathbf{y}(\tau))\,ds\right]\,d\tau\right\|\\
&\leq\int_0^{t_n} \left\| \int_{t_n}^{t_m} K(s-\tau)\mathbf{p}(s-\tau,\mathbf{y}(s),\mathbf{y}(\tau))\,ds \right\|\,d\tau\\
&+\left\|\int_{t_n}^{t_m} \int_\tau^{t_m} K(s-\tau)\mathbf{p}(s-\tau,\mathbf{y}(s),\mathbf{y}(\tau))\,ds\,d\tau\right\|+\left\|\int_{t_n}^{t_m} \mathbf{q}(\mathbf{y}(\tau))\,d\tau\right\|\\
&\leq \int_0^{t_n}  \int_{t_n}^{t_m}\left\| K(s-\tau)\mathbf{p}(s-\tau,\mathbf{y}(s),\mathbf{y}(\tau))\right\|\,ds \,d\tau\\
&+\int_{t_n}^{t_m} \int_\tau^{t_m} \left\|K(s-\tau)\mathbf{p}(s-\tau,\mathbf{y}(s),\mathbf{y}(\tau))\right\|\,ds\,d\tau+\int_{t_n}^{t_m} \left\|\mathbf{q}(\mathbf{y}(\tau))\right\|\,d\tau.
\end{align*}
In the last inequality the third integral tends to zero as $n,m \rightarrow +\infty$, because $\mathbf{q}$ is bounded on $\bar{U}$ and $t_m,t_n \rightarrow T$. Also the first two integrals go to zero as $n,m \rightarrow +\infty$, because of \eqref{cont-Cp2}. The proof is completed.
\end{pf}
\\Now, from Proposition \ref{princ-prolung}, we obtain the following
\begin{Theorem}
Equation \eqref{eq-integrale2} has a unique global solution for all $t\geq 0$.
\end{Theorem}
\begin{pf} First, equations \eqref{qi}--\eqref{qi_end} imply
\begin{align*}
	|q_{1,i1}|\leq 2\beta ||\mathbf{y}||,&\quad |q_{1,i2}|\leq 2\beta ||\mathbf{y}||,\\
	|q_{2,i1}|\leq ||\mathbf{y}||,&\quad |q_{2,i2}|\leq ||\mathbf{y}||,\\
\end{align*}
so that
\begin{align*}
	||\mathbf{q}||\leq \sqrt{2N(4\beta^2+1)}||\mathbf{y}||.
\end{align*}
Then \eqref{eq-integrale1}, \eqref{cont-Cp} yield
\begin{align*}
	||\mathbf{y}||&\leq ||\mathbf{y}_0|| +\sqrt{2N(4\beta^2+1)}\int_0^t||\mathbf{y}(\tau)||\,d\tau\\
	&+\sqrt{2N^3} \left| \left| \varphi\right| \right| \left| \left| \varphi^{\prime}\right| \right|   \int_0^t\int_\tau^t \frac{\left| K(s-\tau)\right|}{\sqrt{s-\tau}}\,ds\,d\tau\\
&=  ||\mathbf{y}_0||+ \sqrt{2N(4\beta^2+1)}\int_0^t||\mathbf{y}(\tau)||\,d\tau +\sqrt{2N^3} \left| \left| \varphi\right| \right| \left| \left| \varphi^{\prime}\right| \right| \int_0^t\int_0^{t-\tau}\frac{\left| K(z)\right|}{\sqrt{z}}\,dz\,d\tau		\\
&\leq ||\mathbf{y}_0||+ \sqrt{2N(4\beta^2+1)}\int_0^t||\mathbf{y}(\tau)||\,d\tau +\sqrt{2N^3} \left| \left| \varphi\right| \right| \left| \left| \varphi^{\prime}\right| \right|\int_0^t\int_0^{+\infty}\frac{\left| K(z)\right|}{\sqrt{z}}\,dz\,d\tau\\
&= ||\mathbf{y}_0||+ \sqrt{2N(4\beta^2+1)}\int_0^t||\mathbf{y}(\tau)||\,d\tau +M\sqrt{2N^3} \left| \left| \varphi\right| \right| \left| \left| \varphi^{\prime}\right| \right|  t,
\end{align*}
where we have set $z=s-\tau$, and $M$ is given by \eqref{int-infinito-C}.
\par Now, for each $0\leq t <T$,
\begin{align*}
	||\mathbf{y}||\leq \left(||\mathbf{y}_0||+M\sqrt{2N^3} \left| \left| \varphi\right| \right| \left| \left| \varphi^{\prime}\right| \right| T\right)+ \sqrt{2N(4\beta^2+1)}\int_0^t||\mathbf{y}(\tau)||\,d\tau,
\end{align*}
so that
\begin{align*}
		||\mathbf{y}||\leq \left(||\mathbf{y}_0||+M\sqrt{2N^3} \left| \left| \varphi\right| \right| \left| \left| \varphi^{\prime}\right| \right|  T\right) e^{t\sqrt{2N(4\beta^2+1)}}
\end{align*}
by the Gronwall's inequality. Since the solution remains bounded, for Proposition \ref{princ-prolung}, it can be continued to all $[0,+\infty)$.
\end{pf}
\section{Asymptotic properties on the linearised system}\label{sec-stab}

In this section we prove some asymptotic properties on the linearised form of system \eqref{sys-cuck-chemo}. To simplify some computations and the following numerical simulations, here we consider the case in which the source term in the chemical signal equation is given by a characteristic function on a ball of radius $R$ centered on each particle, as in \citet{dicostanzo, doi:10.1093/imammb/dqw022}.
Namely, we assume 
\begin{align}\label{function-g}
	g_\mathbf{X}=\xi \displaystyle \sum_{j=1}^{N} \chi_{\mathbf{B}(\mathbf{X}_j,R)},\quad \text{with} \quad \xi>0, 
\end{align}
and 
\begin{align}
\chi_{\mathbf{B}(\mathbf{X}_j,R)}:=
\left\{
\begin{array}{lll}
1,&\text{if}&\mathbf{x}\in\mathbf{B}(\mathbf{X}_j,R):=\left\{\mathbf{x}:\left|\left|\mathbf{x}-\mathbf{X}_j\right|\right|\leq R\right\};\\
0,&&\text{otherwise}.
\end{array}
\right.
\end{align}
With similar computations to those of the previous section, from \eqref{soluz-fond}, \eqref{function-g}, we get
%
\begin{align}\label{calore2D_old}
	f(x_1,x_2,t)= \xi\sum_{j=1}^N\int_0^t \iint_{\mathbf{B}(\mathbf{X}_j(\tau),R)}\frac{e^{-\eta (t-\tau)}e^{-\frac{(x_1-\bar{x}_1)^2+(x_2-\bar{x}_2)^2}{4(t-\tau)D}}}{4\pi(t-\tau)D}\,d\bar{x}_1\,d\bar{x}_2\,d\tau.
\end{align}
%
Performing the change of variables $\tilde{x}_2=\bar{x}_2-X_{j2}(\tau)$, we express the chemotactic gradient $\nabla f=(\partial_{x_1} f,\partial_{x_2} f)$, obtaining
\begin{align}\label{grad_1}
\partial_{x_1} f(x_1,x_2,t) &= -\xi\sum_{j=1}^N\int_0^t\frac{e^{-\eta (t-\tau)}}{4\pi (t-\tau)D}\int_{-R}^{+R}  e^{-\frac{(x_2-X_{j2}(\tau)-\tilde{x}_2)^2}{4(t-\tau)D}}\left(e^{-\frac{\left(x_1-X_{j1}(\tau)-\sqrt{R^2-\tilde{x}_2^2}\right)^2}{4(t-\tau)D} }\right.\nonumber\\
&\left.	-e^{-\frac{\left(x_1-X_{j1}(\tau)+\sqrt{R^2-\tilde{x}_2^2}\right)^2}{4(t-\tau)D} }\right)d\tilde{x}_2\,d \tau,
\end{align}
\noindent
and similarly we can proceed for $\partial_{x_2} f(x_1,x_2,t)$.\\
Hence, substituting  \eqref{grad_1} into \eqref{sys-cuck-chemo}$_1$ we can summarize, for $i=1,\dots,N$, the following system:
\begin{align}\label{Nparticles2D_old}
\left\{
\begin{array}{ll}
	\dot{V}_{i1}&=\displaystyle\frac{\beta}{N}\sum_{j=1}^N\frac{1}{\left(1+\frac{||\mathbf{X}_i-\mathbf{X}_j||^2}{R^2}\right)^{\sigma}}(V_{j1}-V_{i1})\\
	&-\displaystyle\int_0^t C(t-\tau)\sum_{j=1}^N\int_{-R}^{+R}  e^{-\frac{(X_{i2}(t)-X_{j2}(\tau)-\tilde{x}_2)^2}{4(t-\tau)D}}\left(e^{-\frac{\left(X_{i1}(t)-X_{j1}(\tau)-\sqrt{R^2-\tilde{x}_2^2}\right)^2}{4(t-\tau)D} }\right.\\
	&\left.-e^{-\frac{\left(X_{i1}(t)-X_{j1}(\tau)+\sqrt{R^2-\tilde{x}_2^2}\right)^2}{4(t-\tau)D} }\right)d\tilde{x}_2\,d \tau ,\\
	\dot{V}_{i2}&=\displaystyle\frac{\beta}{N}\sum_{j=1}^N\frac{1}{\left(1+\frac{||\mathbf{X}_i-\mathbf{X}_j||^2}{R^2}\right)^{\sigma}}(V_{j2}-V_{i2})\\
	&-\displaystyle\int_0^t C(t-\tau)\sum_{j=1}^N\int_{-R}^{+R}  e^{-\frac{(X_{i1}(t)-X_{j1}(\tau)-\tilde{x}_1)^2}{4(t-\tau)D}}\left(e^{-\frac{\left(X_{i2}(t)-X_{j2}(\tau)-\sqrt{R^2-\tilde{x}_1^2}\right)^2}{4(t-\tau)D} }\right.\\
	&\left.-e^{-\frac{\left(X_{i2}(t)-X_{j2}(\tau)+\sqrt{R^2-\tilde{x}_1^2}\right)^2}{4(t-\tau)D} }\right)d\tilde{x}_1\,d \tau, \\
	\dot{X}_{i1}&=V_{i1},\\
	\dot{X}_{i2}&=V_{i2},
	\end{array}
	\right.
\end{align}

with
       \begin{align}\label{c-stab}
C(t-\tau):=\frac{\gamma\xi e^{-\eta (t-\tau)}}{4\pi (t-\tau)D}.
	\end{align}

We are interested in the equilibrium points that satisfy the condition:
\begin{align}\label{eq-point}
\left\{
\begin{array}{lll}
	\mathbf{X}_i(t)=\mathbf{X}_\text{eq}(t),& \forall i,& \forall t;\\
	\mathbf{V}_i(t)=\mathbf{0}, &\forall i; &
	\end{array}
	\right.
	\quad
	\Leftrightarrow
	\quad
	\mathbf{X}_i(t)=\mathbf{X}_\text{eq}=\text{constant},\quad \forall i.
\end{align}
Equation \eqref{eq-point} means that all particles are in a same position for all times. Now, to make a first-order approximation of \eqref{Nparticles2D_old}, we consider the following Taylor expansions around points \eqref{eq-point}:
\begin{align}\label{lin-fun}
F_1\left(\mathbf{X}_j-\mathbf{X}_i,V_{j1}-V_{i1}\right):=&\frac{1}{\left(1+\frac{\left\|\mathbf{X}_j-\mathbf{X}_i\right\|^2}{R^2}\right)^{\sigma}}(V_{j1}-V_{i1})\\
=&V_{j1}-V_{1i}+\rho_1(\mathbf{X}_j-\mathbf{X}_i,V_{j1}-V_{i1}),\\
F_2\left(t-\tau, X_{i2}(t)-X_{j2}(\tau),\tilde{x}_2\right):=&e^{-\frac{(X_{i2}(t)-X_{j2}(\tau)-\tilde{x}_2)^2}{4(t-\tau)D}}\\
=&e^{-\frac{\tilde{x}_2^2}{4(t-\tau)D}}+\rho_2(t-\tau, X_{i2}(t)-X_{j2}(\tau),\tilde{x}_2),\label{lin-fun1}\\
F_3\left(t-\tau, X_{i1}(t)-X_{j1}(\tau),\tilde{x}_2\right):=&	e^{-\frac{\left(X_{i1}(t)-X_{j1}(\tau)\pm \sqrt{R^2-\tilde{x}_2^2}\right)^2}{4(t-\tau)D}}\\
=&e^{-\frac{R^2-\tilde{x}_2^2}{4(t-\tau)D}} \mp e^{-\frac{R^2-\tilde{x}_2^2}{4(t-\tau)D}}\frac{\sqrt{R^2-\tilde{x}_2^2}}{2(t-\tau)D}\left(X_{i1}(t)-X_{j1}(\tau) \right)\nonumber\\
	+& \rho_3(t-\tau, X_{i1}(t)-X_{j1}(\tau),\tilde{x}_2),\label{lin-fun2}
	\end{align}
where the functions $\rho_1$ and $\rho_3$ contain the nonlinear terms, while $\rho_2$ contains the linear and the nonlinear terms. Similarly we can treat equation \eqref{Nparticles2D_old}$_2$.
\par From \eqref{lin-fun}--\eqref{lin-fun2}, we linearise equation \eqref{Nparticles2D_old}$_1$ in the form
%
\begin{align*}
	\dot{V_{i1}}&=\frac{\beta}{N}\sum_{j=1}^N(V_{j1}(t)-V_{i1}(t))-\displaystyle\int_0^t \frac{C(t-\tau)}{2(t-\tau)D}\int_{-R}^{+R}e^{-\frac{\tilde{x}_2^2}{4(t-\tau)D}}e^{-\frac{R^2-\tilde{x}_2^2}{4(t-\tau)D}}\sqrt{R^2-\tilde{x}_2}\\
	&\cdot\sum_{j=1}^N\left( X_{i1}(t)-X_{j1}(\tau)+X_{i1}(t)-X_{j1}(\tau)\right)\,d\tilde{x}_2\,d\tau\\
	&=\frac{\beta}{N}\sum_{j=1}^N(V_{j1}(t)-V_{i1}(t))-\displaystyle\int_0^t \frac{C(t-\tau) e^{-\frac{R^2}{4(t-\tau)D}}}{(t-\tau)D}\int_{-R}^{+R}\sqrt{R^2-\tilde{x}_2^2}\,d\tilde{x}_2\;\\
	&\cdot\sum_{j=1}^N\left( X_{i1}(t)-X_{j1}(\tau)\right)\,d\tau\\
	&=\frac{\beta}{N}\sum_{j=1}^N(V_{j1}(t)-V_{i1}(t))-\displaystyle\int_0^t \bar{C}(t-\tau)\sum_{j=1}^N\left( X_{i1}(t)-X_{j1}(\tau)\right)\,d\tau,\\
\end{align*}
with
\begin{align}\label{c-barra}
	\bar{C}(t-\tau):=\frac{\pi R^2 C(t-\tau) e^{-\frac{R^2}{4(t-\tau)D}}}{2(t-\tau)D}.
\end{align}
Similarly it can be done for $\dot{V}_{i2}$. Finally, we obtain
\begin{align}\label{sys-lin}
	\left\{
	\begin{array}{ll}
	\dot{V_{i1}}&=\displaystyle\frac{\beta}{N}\displaystyle\sum_{j=1}^N(V_{j1}(t)-V_{i1}(t))-\displaystyle\int_0^t \bar{C}(t-\tau)\sum_{j=1}^N\left( X_{i1}(t)-X_{j1}(\tau)\right)\,d\tau,\\
	\dot{V_{i2}}&=\displaystyle\frac{\beta}{N}\displaystyle\sum_{j=1}^N(V_{j2}(t)-V_{i2}(t))-\displaystyle\int_0^t \bar{C}(t-\tau)\sum_{j=1}^N\left( X_{i2}(t)-X_{j2}(\tau)\right)\,d\tau,\\
	\dot{X_{i1}}&=V_{i1},\\
	\dot{X_{i2}}&=V_{i2}.
	\end{array}
	\right.
\end{align}
We are interested to establish the following time-asymptotic convergence property of \eqref{sys-lin}:
\begin{Theorem}\label{teo-convergence}
Let 
\begin{align}
	\mathbf{X}_\text{CM}(t)&:=\frac{1}{N}\sum_{i=1}^N \mathbf{X}_i(t),\label{cm1}\\
	\mathbf{V}_\text{CM}(t)&:=\frac{1}{N}\sum_{i=1}^N \mathbf{V}_i(t),\label{cm2}
\end{align}
position and velocity of the centre of mass of the system of $N$ particles with same mass, for system \eqref{sys-lin}, for $t\rightarrow +\infty$, the following properties hold:
	\begin{itemize}
		\item[i)] all particles converge towards a same velocity and a same position, that is their centre of mass, as in \eqref{cm1};
		\item[ii)] the velocity of the centre of mass, as in \eqref{cm2}, tends to zero. 
	\end{itemize}
\end{Theorem}
\begin{Remark}
Theorem \ref{teo-convergence} ensures a condition of time-asymptotic flocking such as stated in Definition \ref{def-flocking}. Moreover, we have also the stronger condition that all particles converge asymptotically to their centre of mass and the velocity of the centre of mass decays to zero.
\end{Remark}
Next propositions and lemmas will lead to the proof of Theorem \ref{teo-convergence} at the end of the section. First, it is convenient to introduce the centre of mass system, in which equations \eqref{sys-lin} become a nonautonomous system of ordinary differential equation, decoupled with respect to the $i$-th particle and with respect to the two components of each position and velocity vector. Then the equation of the centre of mass can be studied apart.
\par Starting from \ref{cm1}--\ref{cm2} we define the new variables
\begin{align}
	\bar{\mathbf{X}}_i&:=\mathbf{X}_i-\mathbf{X}_\text{CM},\label{variab-cm}\\
	\bar{\mathbf{V}}_i&:=\mathbf{V}_i-\mathbf{V}_\text{CM}. \label{variab-cm2}
\end{align}
In variables \eqref{variab-cm}--\eqref{variab-cm2} the equilibrium condition \eqref{eq-point} becomes
\begin{align}\label{eq-00}
	(\bar{\mathbf{X}}_i,\bar{\mathbf{V}}_i)=\left(\mathbf{0},\mathbf{0}\right),\quad \forall i=1,\dots,N,
\end{align}
moreover the following identities hold:
\begin{align}
	\sum_{i=1}^N\bar{\mathbf{X}}_i&=\mathbf{0},\label{id-cm1}\\
	\sum_{i=1}^N\bar{\mathbf{V}}_i&=\mathbf{0}.\label{id-cm2}
\end{align}
If $\mathbf{V}_\text{CM}=\left(V_{\text{CM}1},V_{\text{CM}2}\right)$ and $\mathbf{X}_\text{CM}=\left(X_{\text{CM}1},X_{\text{CM}2}\right)$, from \eqref{sys-lin} we have
\begin{align}
	\dot{V}_{\text{CM}1}&=\frac{1}{N}\sum_{i=1}^N \dot{V}_{i1}=\frac{1}{N}\frac{\beta}{N}\sum_{i=1}^N\sum_{j=1}^N(V_{j1}(t)-V_{i1}(t))\nonumber\\
	&-\frac{1}{N}\displaystyle\int_0^t \bar{C}(t-\tau)\sum_{i=1}^N\sum_{j=1}^N\left( X_{i1}(t)- X_{j1}(\tau)\right)\,d\tau\nonumber\\
	&=\frac{1}{N}\frac{\beta}{N}\sum_{i=1}^N(NV_{\text{CM}1}(t)-N V_{i1}(t))\nonumber\\
	&-\frac{1}{N}\displaystyle\int_0^t \bar{C}(t-\tau)\sum_{i=1}^N\left(N X_{i1}(t)- NX_{\text{CM}1}(\tau)\right)\,d\tau\nonumber\\
	&=\frac{1}{N}\frac{\beta}{N}(N^2V_{\text{CM}1}(t)-N^2 V_{\text{CM}1}(t))\nonumber\\
	&-\frac{1}{N}\displaystyle\int_0^t \bar{C}(t-\tau)\left(N^2 X_{\text{CM}1}(t)- N^2X_{\text{CM}1}(\tau)\right)\,d\tau\nonumber\\
	&=-\displaystyle\int_0^t \bar{C}(t-\tau)\left(N X_{\text{CM}1}(t)- NX_{\text{CM}1}(\tau)\right)\,d\tau,\label{V-punto-CM}
\end{align}
where we have used definitions \eqref{cm1}--\eqref{cm2}.
The same holds for $\dot{V}_{\text{CM}2}$:
\begin{align}
	\dot{V}_{\text{CM}2}=-\displaystyle\int_0^t \bar{C}(t-\tau)\left(N X_{\text{CM}2}(t)- NX_{\text{CM}2}(\tau)\right)\,d\tau.\label{V-punto-CM2}
\end{align}
In the variables $(\bar{\mathbf{X}}_i,\bar{\mathbf{V}}_i)$, taking into account \eqref{V-punto-CM}--\eqref{V-punto-CM2}, equations \eqref{sys-lin}$_{1,3}$ become
\begin{align*}
\dot{\bar{V}}_{i1}&=-\dot{V}_{\text{CM}1}+\frac{\beta}{N}\sum_{j=1}^N(\bar{V}_{j1}(t)-\bar{V}_{i1}(t))\\
	&-\displaystyle\int_0^t \bar{C}(t-\tau)\sum_{j=1}^N\left( \bar{X}_{i1}(t)-\bar{X}_{j1}(\tau)+X_{\text{CM}1}(t)-X_{\text{CM}1}(\tau)\right)\,d\tau\\
	&=\frac{\beta}{N}\left(\sum_{j=1}^N\bar{V}_{j1}(t)-\sum_{j=1}^N\bar{V}_{i1}(t)\right)\\
	&-\displaystyle\int_0^t \bar{C}(t-\tau)( N\bar{X}_{i1}(t)-\sum_{j=1}^N\bar{X}_{j1}(\tau)+NX_{\text{CM}1}(t)-NX_{\text{CM}1}(\tau) \\
	& -NX_{\text{CM}1}(t)+NX_{\text{CM}1}(\tau))\,d\tau\\
	&=-\beta \bar{V}_{i1}(t)- N\left(\displaystyle\int_0^t \bar{C}(t-\tau)\,d\tau\right)\bar{X}_{i1}(t),
	\end{align*}
where we have used \eqref{id-cm1}--\eqref{id-cm2}, and
\begin{align*}
\dot{\bar{X}}_{i1}&=-\dot{X}_{\text{CM}1}+\bar{V}_{i1}+V_{\text{CM}1}=\bar{V}_{i1}.
\end{align*}
Similarly for $\dot{\bar{V}}_{i2}$ and $\dot{\bar{X}}_{i2}$. \par Finally, we can write the following system:
\begin{align}\label{sys-lin-CM}
\left\{
	\begin{array}{l}
	\dot{\bar{V}}_{i1}=-\beta\bar{V}_{i1}(t)-g(t)\bar{X}_{i1}(t),\\
	\dot{\bar{V}}_{i2}=-\beta\bar{V}_{i2}(t)-g(t)\bar{X}_{i2}(t),\\
	\dot{\bar{X}}_{i1}=\bar{V}_{i1},\\
	\dot{\bar{X}}_{i2}=\bar{V}_{i2},
	\end{array}
	\right.
\end{align}
where
\begin{align*}
	g(t):=N\int_0^t \bar{C}(t-\tau)\,d\tau.
\end{align*}
Now we will prove the uniform asymptotic stability of equilibrium \eqref{eq-00} providing a suitable Lyapunov function for system \eqref{sys-lin-CM}. For simplicity, system \eqref{sys-lin-CM} can be written, for each particle and for each component, as a planar system in the variable $\mathbf{y}=\left(V,X\right)$:
\begin{align}\label{sys-lin-CM-VX}
\left\{
	\begin{array}{l}
	\dot{V}=-\beta V-g(t)X,\\
	\dot{X}=V,
	\end{array}
	\right.
\end{align}
with
\begin{align}\label{g-t}
	g(t):=N\int_0^t \bar{C}(t-\tau)\,d\tau.
\end{align}
In relation to \eqref{sys-lin-CM-VX} we prove the following two propositions.
\begin{Proposition}\label{prop-lya}
Fixed a $\bar{t}>0$, the system \eqref{sys-lin-CM-VX}, admits a Lyapunov function $U(t,\mathbf{y})$ with the properties:
\begin{enumerate}[a)]
\item
\begin{align}\label{lya1}
	k_2 \left\|\mathbf{y}\right\|^2 \leq U(t,\mathbf{y})\leq k_1 \left\|\mathbf{y}\right\|^2;
\end{align}
\item
\begin{align}\label{lya3}
	\dot{U}(t,\mathbf{y})=\frac{\partial U}{\partial t}+\frac{\partial U}{\partial V}\dot{V}+\frac{\partial U}{\partial X}\dot{X} \leq -k_3 \left\|\mathbf{y}\right\|^2;
\end{align}
\end{enumerate}
for all $t\geq \bar{t}$, where $k_1$, $k_2$, and $k_3$ are positive constants.
\end{Proposition}
\begin{pf} Let $\bar{t}>0$, we define the Lyapunov function
\begin{align*}
U(t,V,X):=(V^2+kXV+g(t)X^2)\psi(t),
\end{align*}
where
\begin{align}
\psi(t)&:=e^{-\frac{g(t)}{\underline{g}}},\label{psi-t}\\
\underline{g}&:=\inf_{t\geq \bar{t}} g(t)=N\int_0^{\bar{t}} \bar{C}(t-\tau)\,d\tau,\quad \overline{g}:=\sup_{t\geq \bar{t}} g(t)=N\int_0^{+\infty} \bar{C}(t-\tau)\,d\tau,\\
\underline{\psi}&:=\inf_{t\geq \bar{t}} \psi=e^{-\frac{\bar{g}}{\underline{g}}},\quad \overline{\psi}:=\sup_{t\geq \bar{t}} \psi=e^{-1},\\
\overline{\overline{\psi}}&:=\sup_{t\geq \bar{t}} |\dot{\psi}|=\sup_{t\geq \bar{t}} \frac{e^{-\frac{g(t)}{\underline{g}}}\dot{g}}{\underline{g}},\\
k&:=\min\left[\frac{\;\underline{\psi}\;}{\overline{\psi}},\frac{\underline{g}\,\underline{\psi}}{\overline{\psi}},
\frac{2\beta \underline{g}\,\underline{\psi}^2}{2 \underline{g}\,\underline{\psi}\overline{\psi}+(\overline{\overline{\psi}}+\beta\overline{\psi})^2}\label{kappa-min}
\right].
\end{align}
Similar functional can be found in \citet{onitsuka2011}. In the following equations we consider the inequalities:
\begin{align*}
	-\frac{X^2+V^2}{2} \leq XV\leq \frac{X^2+V^2}{2}.
\end{align*}
Then, because of $g(t)$ is an increasing function, $\psi(t)$ is nonincreasing, so
\begin{align}\label{psi-punto}
	\dot{\psi}\leq 0,
\end{align}
and finally
\begin{align}
\dot{\psi}g+\psi\dot{g}&=-e^{-\frac{g(t)}{\underline{g}}}\frac{\dot{g}}{\underline{g}}g+e^{-\frac{g(t)}{\underline{g}}}\dot{g}\nonumber\\
	&=\dot{g}e^{-\frac{g(t)}{\underline{g}}}\left(1-\frac{g(t)}{\underline{g}}\right)\nonumber\\
	&\leq 0, \quad \forall\, t\geq \bar{t}.\label{der-psi-g}
\end{align}
To prove the second inequality in a) we consider 
\begin{align*}
	U(t,V,X)&\leq \left(V^2+ k\frac{X^2+V^2}{2}+\overline{g}X^2\right)\overline{\psi}\\
	&=\left[\left(1+\frac{k}{2}\right)V^2+ \left(\overline{g}+\frac{k}{2}\right)X^2\right]\overline{\psi}\\
	&\leq k_1 \left\|\mathbf{y}\right\|^2,
\end{align*}
where
\begin{align*}
	k_1:=\max\left[\left(1+\frac{k}{2}\right)\overline{\psi},\left(\overline{g}+\frac{k}{2}\right) \overline{\psi} \right].
\end{align*}
To prove the first part of a) we consider
\begin{align*}
	U(t,V,X)&\geq \underline{\psi} V^2+\underline{g}\,\underline{\psi}X^2-k\overline{\psi}\frac{X^2+V^2}{2}\\
	&=V^2\left(\underline{\psi}-\frac{k \overline{\psi}}{2}\right)+X^2\left(\underline{g}\,\underline{\psi}-\frac{k \overline{\psi}}{2}\right)\\
	&\geq_{\scriptscriptstyle\eqref{kappa-min}} V^2\frac{\underline{\psi}}{2}+X^2\frac{\underline{g}\,\underline{\psi}}{2}\\
	&\geq k_2\left\|\mathbf{y}\right\|^2,
\end{align*}
where
\begin{align*}
	k_2:=\min\left[\frac{\;\underline{\psi}\;}{2}, \frac{\underline{g}\,\underline{\psi}}{2}\right]
\end{align*}
To prove b) we consider the following inequalities:
\begin{align*}
	\dot{U}(t,V,X)&=\dot{\psi}(V^2+kXV+gX^2)+\psi\left[2V(-\beta V -g X)+kV^2
	+kX (-\beta V -gX)+\dot{g}X^2+2gXV\right]\\
&=\left(\dot{\psi}-2\beta \psi +k\psi\right)V^2+\left(k\dot{\psi}-\beta k\psi\right)XV+\left(\dot{\psi}g+\psi\dot{g}-kg\psi\right) X^2\\
&\leq_{\scriptscriptstyle\eqref{psi-punto},\eqref{der-psi-g}} \left(-2\beta \psi +k\psi\right)V^2+\left(k\dot{\psi}-\beta k\psi\right)XV-kg\psi X^2\\
&\leq \left(-2\beta \psi +k\psi\right)V^2+\left(k|\dot{\psi}|+\beta k\psi\right)|X||V|-kg\psi X^2\\
&\leq \left(-2\beta \underline{\psi} +k\overline{\psi}\right)V^2+\left(k\overline{\overline{\psi}}+\beta k\overline{\psi}\right)|X||V|-k\underline{g}\,\underline{\psi} X^2\\
&=\left(-2\beta \underline{\psi} +k\overline{\psi}\right)V^2-\frac{k\underline{g}\, \underline{\psi}}{2}\left[|X|-\frac{|V|\left(\overline{\overline{\psi}}+\beta \overline{\psi}\right)}{\underline{g}\,\underline{\psi}}\right]^2-\frac{k\underline{g}\,\underline{\psi}}{2} X^2 
+\frac{k\left(\overline{\overline{\psi}}+\beta \overline{\psi}\right)^2}{2 \underline{g}\,\underline{\psi}}V^2\\
&\leq  \left(-2\beta \underline{\psi} +k\overline{\psi}+\frac{k\left(\overline{\overline{\psi}}+\beta \overline{\psi}\right)^2}{2 \underline{g}\,\underline{\psi}} \right)V^2-\frac{k\underline{g}\,\underline{\psi}}{2} X^2\\
&\leq_{\scriptscriptstyle\eqref{kappa-min}} -\left(\beta \underline{\psi} V^2+\frac{k\underline{g}\,\underline{\psi}}{2} X^2\right)\\
&\leq - k_3 \left\|\mathbf{y}\right\|^2,
\end{align*}
where
\begin{align*}
	k_3:=\min \left[\beta \underline{\psi}, \frac{k\underline{g}\,\underline{\psi}}{2} \right].
\end{align*}
This completes the proof.
\end{pf}
\\Starting from Proposition \ref{prop-lya} we can state the following
%
%
\begin{Proposition}\label{teo-as-stab}
The equilibrium point $\left(V,X\right)=(0,0)$ of the linearised system \eqref{sys-lin-CM-VX} is globally uniformly asymptotically stable with exponential rate of convergence.
\end{Proposition}
\begin{pf} Inequalities \eqref{lya3} and \eqref{lya1} imply that $U$ satisfies the differential inequality
\begin{align*}
	\dot{U}\leq -\frac{k_3}{k_1}U,\quad \forall\, t\geq \bar{t}.
\end{align*}
By the Gronwall's inequality,
\begin{align*}
	U(t,\mathbf{y}(t))\leq U\left(\bar{t},\mathbf{y}\left(\bar{t}\right)\right)e^{-(k_3/k_1)(t-\bar{t})}.
\end{align*}
Then, using again \eqref{lya1}, we have
\begin{align}
	\left\|\mathbf{y}(t)\right\|&\leq \left(\frac{U(t,\mathbf{y}(t))}{k_2}\right)^{1/2}\nonumber\\
	&\leq \left(\frac{U\left(\bar{t},\mathbf{y}\left(\bar{t}\right)\right)e^{-(k_3/k_1)(t-\bar{t})}}{k_2}\right)^{1/2}\nonumber\\
	& \leq \left(\frac{k_1 \left\|\mathbf{y}\left(\bar{t}\right)\right\|^2e^{-(k_3/k_1)(t-\bar{t})}}{k_2}\right)^{1/2}\nonumber\\
	 &=\left(\frac{k_1}{k_2}\right)^{1/2} \left\|\mathbf{y}\left(\bar{t}\right)\right\| e^{-(k_3/(2k_1))(t-\bar{t})}.\label{exp-stab}
\end{align}
From \eqref{exp-stab} conditions for the uniform asymptotic stability are satisfied with exponential convergence. The proof is completed.
\end{pf}
\begin{Remark}\label{oss:teo:parte1}
Returning to system \eqref{sys-lin-CM}, Proposition \ref{teo-as-stab} can be applied for each particle and for each component of the position and velocity vectors. Recalling transformations \eqref{variab-cm}--\eqref{variab-cm2}, this proves the first part of Theorem \ref{teo-convergence}.
\end{Remark}
%
%
%
\par Now, to prove the second part of Theorem \ref{teo-convergence}, we investigate equation \eqref{V-punto-CM} for the motion of the centre of mass. Taking into account that $\dot{X}_{\text{CM}1}=V_{\text{CM}1}$, it assumes the form of a Volterra Integro-Differential Equations (VIDEs) of the type
\begin{equation}\label{eq:ezio}
\dot{v}(t)=-\int_0^tc(t-\tau)\int_{\tau}^tv(s)dsd\tau,
\end{equation}
with
\begin{equation}\label{eq:c}
c(u)=c_1\frac{e^{-c_2u}e^{-\frac{c_3}{u}}}{u^2},
\end{equation}
$c_1,\;c_2$ and $c_3$ being positive constants, that arise from \eqref{c-barra} and \eqref{c-stab}. The same holds for equation \eqref{V-punto-CM2}.
\par While existence and uniqueness of the solution of \eqref{eq:ezio} follow from classic techniques, here we want to study its asymptotic properties, which give the asymptotic behaviour of the centre of mass. 
By Dirichlet formula, equation \eqref{eq:ezio} becomes
\[\dot{v}(t)=-\int_0^t\int_0^s c(t-\tau)d\tau v(s)ds,\]
or equivalently
\begin{equation}\label{eq:eq}
\dot{v}(t)=-\int_0^tK(t,s)v(s)ds,\;\;K(t,s)=\int_{t-s}^{t}c(u)du.
\end{equation}
The analytic form of the kernel $K$ is not known. By observing that $\lim_{u\to\infty}u^{2}c(u)=0,$ we can write $K(t,s)=F(t)-\alpha+\alpha-F(t-s),$ where 
\begin{equation}\label{eq:F}
F(t)=\int_0^tc(u)du,
\end{equation}
and
\begin{equation}\label{eq:limF}
\alpha=\alpha(c_1,c_2,c_3)=\int_0^{\infty}c(u)du=\lim_{t\to\infty}F(t).
\end{equation}
With this notation equation \eqref{eq:eq} can be written in the form
(VIDEs) of the type \citet[eq. (9.9)]{lubich83}, i.e.
\begin{equation}\label{eq:Lub}
\dot{v}(t)=\bar{f}(t)+\int_0^tB(t-s)v(s)ds,
\end{equation}
with
\begin{equation}\label{eq:Lubf}
\bar{f}(t)=-(F(t)-\alpha)\int_0^tv(s)ds,
\end{equation}
and
\begin{equation}\label{eq:LubB}
B(t)=F(t)-\alpha.
\end{equation}
Our aim now is to apply the theorem by Miller and Grossman in \citet[Th. 9.2]{lubich83}, that we report here adapted to a scalar equation of the type \eqref{eq:Lub}.
\begin{Theorem}[\textbf{Miller and Grossman}]\label{th:lub}
	Assume that, in equation \eqref{eq:Lub}, the kernel $B$ is in $L^1(0,\infty).$ Then $v(t)\to 0$ whenever $\bar{f}(t)\to 0$ if and only if 
	\begin{equation}\label{eq:tesi}
	w(z)=z-\int_{0}^{+\infty}B(u)e^{-zu}du\neq 0,\;\;\text{Re}(z)\geq 0.
	\end{equation}
\end{Theorem}
The following results are then the fundamental premise for applying Theorem \ref{th:lub}.
\begin{lemma}\label{th:th1}
	The function $B(t)$, defined in \eqref{eq:LubB} with $F(t)$ given in \eqref{eq:F}, satisfies \eqref{eq:tesi}.
\end{lemma}
\begin{pf}
	For $z=0,$ we have $w(0)=\int_0^{+\infty}B(u)du<0,$ since $F(t)<\alpha.$ When $z\neq 0,$ according to the elementary properties of the Laplace transform, we get
	\[w(z)=z-\frac{1}{z}\int_0^{+\infty}c(u)e^{-zu}du +\frac{\alpha}{z}=\frac{1}{z}\varphi(z),\]
	with $\varphi(z)=z^2-\int_0^{+\infty}c(u)e^{-zu}du+\alpha.$
	Then, with $z=x+iy,$  we have
	$\varphi(z)=\varphi(x,y)=a(x,y)+ib(x,y),$ with
	\begin{eqnarray*}
	a(x,y)&=&x^2-y^2+\alpha-\int_0^{+\infty}c(u)e^{-xu}\cos(yu)du\\
	b(x,y)&=&2xy+\int_0^{+\infty}c(u)e^{-xu}\sin(yu)du.
	\end{eqnarray*}
	Consider the case $x>0,$ $y>0.$ Let $p>0$ be an arbitrary constant and denote by $\tilde x=\frac{x}{\sqrt{p}}$ and $\tilde y=\frac{y}{\sqrt{p}}.$ Then 
	\[a(x,y)=0\Leftrightarrow 
	\tilde a(\tilde x,\tilde y)=\tilde x^2-\tilde y^2+\frac{\alpha}{p}-\frac{1}{p}\int_0^{+\infty}c(u)e^{-\sqrt{p}\tilde x u}\cos(\sqrt{p}\tilde y u)du=0,\] and
	$$b(x,y)=0\Leftrightarrow 
	\tilde b(\tilde x,\tilde y)=2\tilde x\tilde y+\frac{1}{p}\int_0^{+\infty}c(u)e^{-\sqrt{p}\tilde x u}\sin(\sqrt{p}\tilde y u)du=0.$$
	Taking into account that
	\[\left|\int_0^{+\infty}c(u)e^{- x u}\cos({ y u})du\right|\leq \alpha,\;\;\forall x>0,\;y\in \mathbb R,\]
	and 
		\[\left|\int_0^{+\infty}c(u)e^{- x u}\sin( y u)du\right|\leq \alpha,\;\;\forall x>0,\;y\in \mathbb R,\] we deduce that 
		\[\tilde x^2-\tilde y^2< \tilde a(\tilde x,\tilde y)<\tilde x^2-\tilde y^2+2\frac{\alpha}{p},\]
		\[\tilde b(\tilde x,\tilde y)>2\tilde x\tilde y-\frac{\alpha}{p},\]
		and hence
		$\tilde a(\tilde x,\tilde y)>0$ for $0<\tilde y<\tilde x,$ and $\tilde a(\tilde x,\tilde y)<0$ for $\tilde y>\sqrt{\tilde x^2+2\frac{\alpha}{p}},$ whereas $\tilde b(\tilde x,\tilde y)>0$ for $\tilde y>\frac{\alpha}{2p\tilde x}.$ This assures that $\tilde a$ and $\tilde b$  may vanish simultaneously only in the white area $\mathcal A$ plotted in Figure \ref{fig:1}, which, choosing $p$ properly, can be made small and included in a  square having side of length 1. Finally, by plotting $\tilde a(\tilde x,\tilde y)$ and $\tilde b(\tilde x,\tilde y)$ for $0< \tilde x<1,$ $0<\tilde y<1,$ with the help of Mathematica, we can verify that it definitely results that they never vanish at the same time. Hence, also $a(x,y)$ and $b(x,y)$ are never both zero at a point $(x,y)$ and \eqref{eq:tesi} is true. \\The case $x\geq 0,$ $y<0$ can be treated analogously, since $a(x,y)=a(x,-y)$ and $b(x,y)=-b(x,-y).$	
\begin{figure}[tbh]
	\begin{center}
		$\begin{array}{c}
		\mbox{\includegraphics[width=0.75\textwidth]{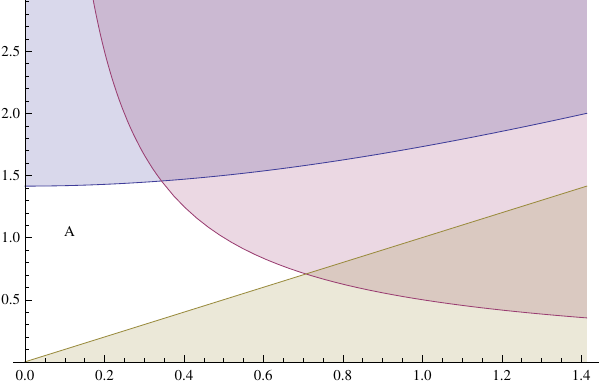}}
		\end{array}$
		\caption{\small    $\mathcal A:$ vanishing area for $\varphi(\tilde x,\tilde y)$.}
		\label{fig:1}
	\end{center}
\end{figure}
\end{pf}

\vspace{1 cm}

\begin{lemma}\label{th:th2}
	Assume that, in equation \eqref{eq:Lub},
	\begin{description}
		\item[i)] $\exists$ $M>0$ and $c_4\leq c_2$ such that $|v(u)|<Me^{c_4t},$
	\end{description}
	where $c_2$ is one of the positive constants appearing in the definition of kernel $B$ by \eqref{eq:c}, \eqref{eq:F} and \eqref{eq:LubB}, then $\lim_{u\to +\infty}|v(u)|=0.$
\end{lemma}
\begin{pf}
	In the assumptions $i),$ the function $\bar{f}$ in \eqref{eq:Lubf} satisfies\[|\bar{f}(t)|\leq (\alpha-F(t))\frac{M}{c_4}(e^{c_4t}-1).\]
	By the de l'Hospital rule,
	\[\lim_{t\to +\infty}(\alpha-F(t))e^{c_4t}=\frac{1}{c_4}\lim_{t\to +\infty}c(t)e^{c_4t}.\]
	Then by recalling the expression \eqref{eq:c} of $c(t),$ the hypothesis $c_4\leq c_2$ and \eqref{eq:limF}, we get
	\begin{equation}\label{eq:fto0}
	\lim_{t\to +\infty}\bar{f}(t)=0.
	\end{equation}
	By analogous considerations it can be shown that $\lim_{t\to +\infty}(F(t)-\alpha)t^2=0,$ so that 
	\begin{equation}\label{eq:Bl1}
	B(t)\in L^1.
	\end{equation}
	Lemma \ref{th:th1}, together with \eqref{eq:fto0} and \eqref{eq:Bl1} assure that all the assumptions of \citet[Th. 9.2]{lubich83} are accomplished, so that $\lim_{t\to +\infty}v(t)=0.$
	\end{pf}
\begin{Remark}
	Lemma \ref{th:th2} requires that the solution of \eqref{eq:eq} is bounded by an increasing exponential function. So, the class of functions involved in this result is quite large and assures that any bounded solution of equation  \eqref{eq:Lub}, and thus of \eqref{eq:ezio}, vanishes at infinity. Nevertheless we want to prove that this is true for any solution of \eqref{eq:ezio}. 
\end{Remark}
The following result, which generalizes Lemma \ref{th:th2}, considerably relaxes the hypotheses on the function $v(t).$
\begin{lemma}\label{th:th3}
	Assume that, in equation \eqref{eq:Lub},
	\begin{equation}\label{eq:M0}
	\exists M_0>0 \mbox{ and } n\in \mathbb N: |v(t)|\leq M_0e^{nc_2t}.
	\end{equation}
	Then $\lim_{t\to +\infty} v(t)=0.$
\end{lemma}
\begin{pf}
	If we multiply both sides of equation \eqref{eq:Lub} by $e^{-(n-i)c_2t}$ and add and subtract the quantity $(n-i)c_2e^{-(n-i)c_2t}v(t),$ we obtain 
	\[\begin{split} e^{-(n-i)c_2t}\dot{v}(t)-(n-i)c_2e^{-(n-i)c_2t}v(t)=e^{-(n-i)c_2t}f(t)-(n-i)c_2e^{-(n-i)c_2t}v(t)\\+\int_0^tB(t-s)e^{-(n-i)c_2(t-s)}e^{-(n-i)c_2s}v(s)ds.\end{split}\]
	Denoting by $v_i(t)=e^{-(n-i)c_2t}v(t),$ $B_i(t)=e^{-(n-i)c_2t}B(t)$ and
	\begin{equation}\label{eq:fi}
	\bar{f}_i(t)=e^{-(n-i)c_2t}\bar{f}(t),
	\end{equation}
	this reads
	\begin{equation}\label{eq:vi}
	\dot{v}_i(t)=\bar{f}_i(t)-(n-i)c_2v_i(t)+\int_0^tB_i(t-s)v_i(s)ds.
	\end{equation}
	Now we want to prove that, for $i=0,\ldots, n-1,$ equation \eqref{eq:vi} satisfies the hypotheses of \citet[Th. 9.2]{lubich83}. First of all we note that $B_i\in L^1[0,+\infty),$ for all $i=0,\ldots, n-1.$ Consider 
	\begin{equation}\label{eq:zi}
	z+(n-i)c_2-\int_0^{+\infty}B_i(u)e^{-zu}du,
	\end{equation}
	which, setting $\zeta=z+(n-i)c_2,$ becomes
	\begin{equation}\label{eq:zeta}
	\zeta-\int_0^{+\infty}B(u)e^{-\zeta u}du.
	\end{equation}
In Lemma \ref{th:th1} ve have proved that the expression in \eqref{eq:zeta} is nonzero for all $\text{Re}(\zeta)\geq 0,$ hence also \eqref{eq:zi} is nonzero for all $\text{Re}(z)\geq 0.$ \\ Now consider $f_i(t).$ For $i=1,$ from \eqref{eq:M0} and \eqref{eq:fi}, \[|\bar{f}_1(t)|\leq (\alpha-F(t))\left(\frac{M_0}{nc_2}e^{nc_2t}-1\right)e^{-(n-1)c_2t}.\] According to the proof of Lemma \ref{th:th2}, it is easy to see that $\lim_{t\to +\infty}(\alpha-F(t))e^{nc_2t}e^{-(n-1)c_2t}=\lim_{t\to +\infty}(\alpha-F(t))e^{c_2t}=0.$ Then, $\lim_{t\to +\infty}\bar{f}_1(t)=0$ and, in view of \citet[Th. 9.2]{lubich83}, we have $\lim_{t\to +\infty}v_1(t)=0.$ This in turn implies that
	\begin{equation}\label{eq:M1}
	\exists M_1>0 : |v(t)|\leq M_1e^{(n-1)c_2t}.
	\end{equation}
	Now consider $i=2,$  and assume that \eqref{eq:M1} holds. By analogous considerations on $\bar{f}_2(t)$ we obtain that $\exists M_2>0 : |v(t)|\leq M_2e^{(n-2)c_2t}.$ Proceeding further for $i=3,\;4,\ldots,n-1,$ we obtain that $\exists M_{n-1}>0 : |v(t)|\leq M_{n-1}e^{c_2t},$ thus  satisfying the hypotheses of Lemma \ref{th:th2} with $M=M_{n-1}$ and $c_4=c_2.$ 
	\end{pf}
\begin{Proposition}\label{th:th4}
	For any triple of positive constants $c_1,\;c_2$ and $c_3,$ the solution $v(t)$ of equation   \eqref{eq:ezio} satisfies $\lim_{t\to +\infty}v(t)=0.$
\end{Proposition}
\begin{pf} 
	By integrating both sides of the VIDE in \eqref{eq:eq} we obtain the Volterra integral equation
	\[v(t)=v(0)-\int_0^tA(t,s)v(s)ds,\]
	with $A(t,s)=\int_s^t\int_{\tau-s}^{\tau}c(u)dud\tau.$
	Observe that, since $\frac{e^{-\frac{c_3}{u}}}{u^2}\leq \frac{4}{c_3^2e^2},$ it holds that $c(u)\leq\frac{4c_1}{c_3^2e^2}e^{-c_2u}.$ Therefore, for $\tilde v(t)= e^{-c_2t}v(t),$ the following inequality holds.
	\[|\tilde v(t)|\leq e^{-c_2t}|v(0)|+\int_0^te^{-c_2(t-s)}\int_s^t\int_{\tau-s}^{\tau}\frac{4c_1}{c_3^2e^2}e^{-c_2u}dud\tau|\tilde v(s)|ds,\]
	hence,
	\[|\tilde v(t)|\leq |v(0)|+\frac{4c_1}{c_2^2c_3^2e^2}\int_0^t\Lambda(t,s)|\tilde v(s)|ds,\]
	where $\Lambda(t,s)=e^{-2c_2t}(e^{c_2s}-1)(e^{c_2t}-e^{c_2s})\leq 1.$  So, by Gronwall inequality (see e.g. \citealp[pg. 79]{brubook}) 
	\[| v(t)|\leq |v(0)|e^{\left(\frac{4c_1}{c_2^2c_3^2e^2}+c_2\right)t},\]
	and the result follows from Lemma \ref{th:th3}.
	\end{pf}
\begin{Remark}\label{oss:teo:parte2}
Proposition \ref{th:th4} holds for both components $V_{\text{CM}1}$, $V_{\text{CM}2}$ with equations \eqref{V-punto-CM}, \eqref{V-punto-CM2}
\end{Remark}

\begin{pf}\textbf{(Theorem \ref{teo-convergence})} First part of the theorem is proved by Proposition \ref{teo-as-stab} and Remark \ref{oss:teo:parte1}, second part follows from Proposition \ref{th:th4} and Remark \ref{oss:teo:parte2}.
\end{pf}

\section{Numerical simulations}\label{sec:numerical-chemo}
In this section we present one of the performed numerical simulations to show the dynamical behaviour of the model introduced in Section \ref{sec:basic-flock-chemo}. Numerical results are compared with the analytical ones presented in Section \ref{sec:local}. \\Numerical tests are performed and shown in nondimensional form using the following dimensionless quantities:
\begin{equation*}\label{sys-numerica}
\begin{split}
  t^*:=t \eta,\quad \mathbf{X}^*:=\frac{\mathbf{X}}{R},\quad &f^*:=\frac{f}{f_{\max}},\quad \beta^*:=\frac{\beta}{\eta},
	\quad\gamma^*:=\frac{\gamma f_{\max}}{R^2\eta^2},\quad D^*:=\frac{D }{R^2\eta},\quad \xi^*:=\frac{\xi}{f_{\max}\eta},
\end{split}
\end{equation*}
where $f_{\max}$ is the maximum concentration of signal $f$. With these definitions, system \eqref{sys-cuck-chemo} can be written as
\begin{align}\label{sys-cuck-chemo-adim}
\left\{
	\begin{array}{l}
	\dot{\mathbf{V}}_i=\displaystyle\frac{\beta}{N}\sum_{j=1}^N\frac{1}{\left(1+\left\|\mathbf{X}_i-\mathbf{X}_j\right\|^2\right)^{\sigma}}(\mathbf{V}_j-\mathbf{V}_i)+\gamma \nabla f (\mathbf{X}_i),\\
	\dot{\mathbf{X}}_i=\mathbf{V}_i,\\
	\partial_t f=D\Delta f+\xi \displaystyle\sum_{j=1}^N \chi_{\mathbf{B}\left(\mathbf{X}_j,1\right)}-f,
	\end{array}
	\right.
\end{align}
where we have dropped, for simplicity, the asterisks for the nondimensional quantities. Notice that, due to the choice of $R$ as characteristic length, the dimensionless particle ray turns out to be a unit value.
\subsection{Numerical methods}
The numerical approximation scheme used here employs a 2D finite difference method on a spatial domain $\Omega:=[a,b]\times[c,d]$ with periodic boundary conditions.
\par For the parabolic equation \eqref{sys-cuck-chemo-adim}$_3$, in order to eliminate the stiff term $-f$ we perform the classic exponential transformation, and we apply a central difference scheme in space and the parabolic Crank-Nicolson scheme in time.
\par For equations \eqref{sys-cuck-chemo-adim}$_{1}$ we adopt a one step IMEX method, putting in implicit the term depending on the velocities and in explicit the gradient term \citep{hundsdorfer}.
%
%
%
%
\subsection{Numerical test}\label{sub-sec-tests}
\begin{figure}[!ht]
\centering
\subfigure[]{\includegraphics[width=0.45\textwidth]{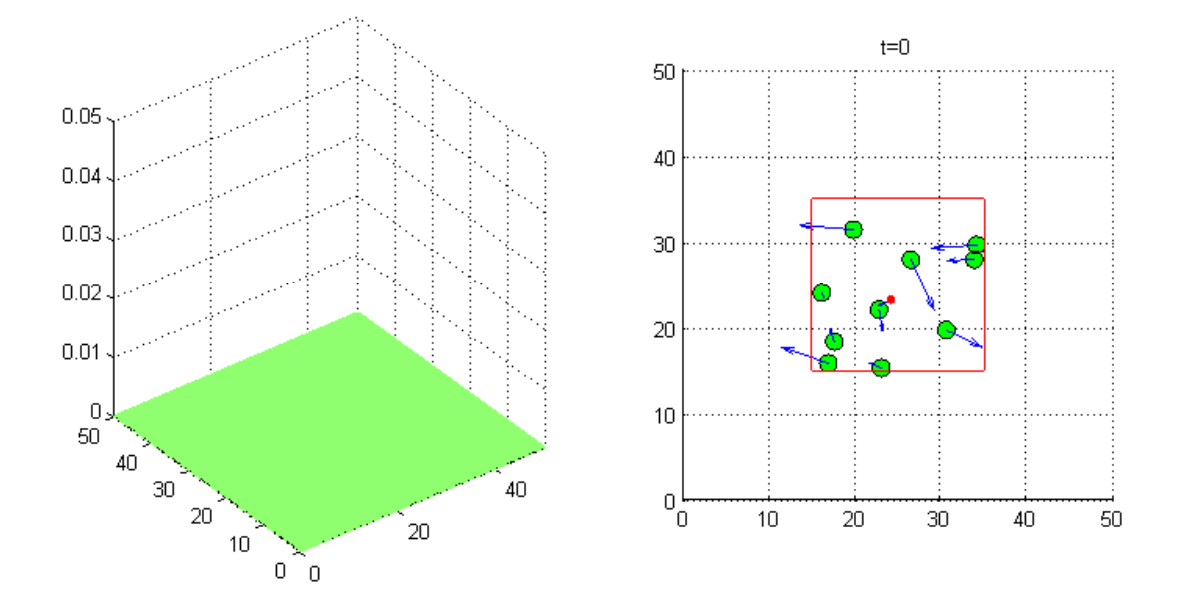}}
\subfigure[]{\includegraphics[width=0.45\textwidth]{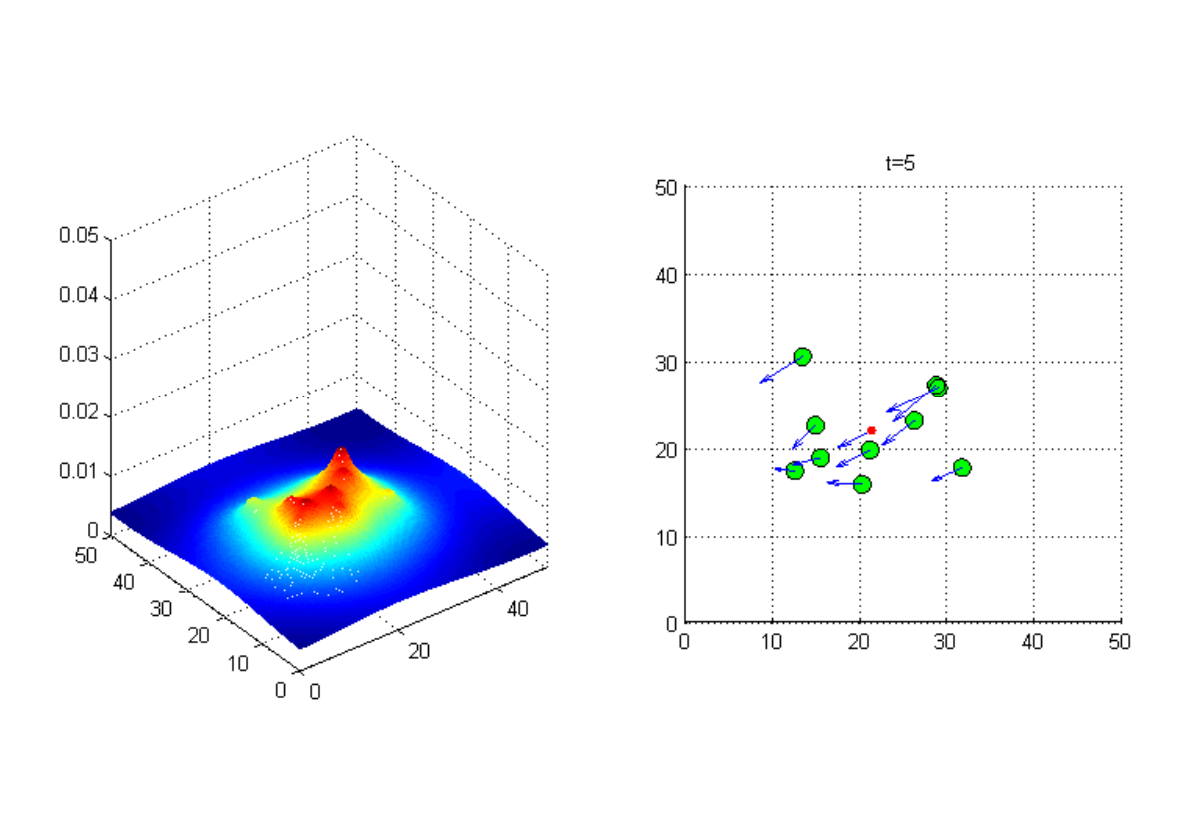}}\\
\subfigure[]{\includegraphics[width=0.45\textwidth]{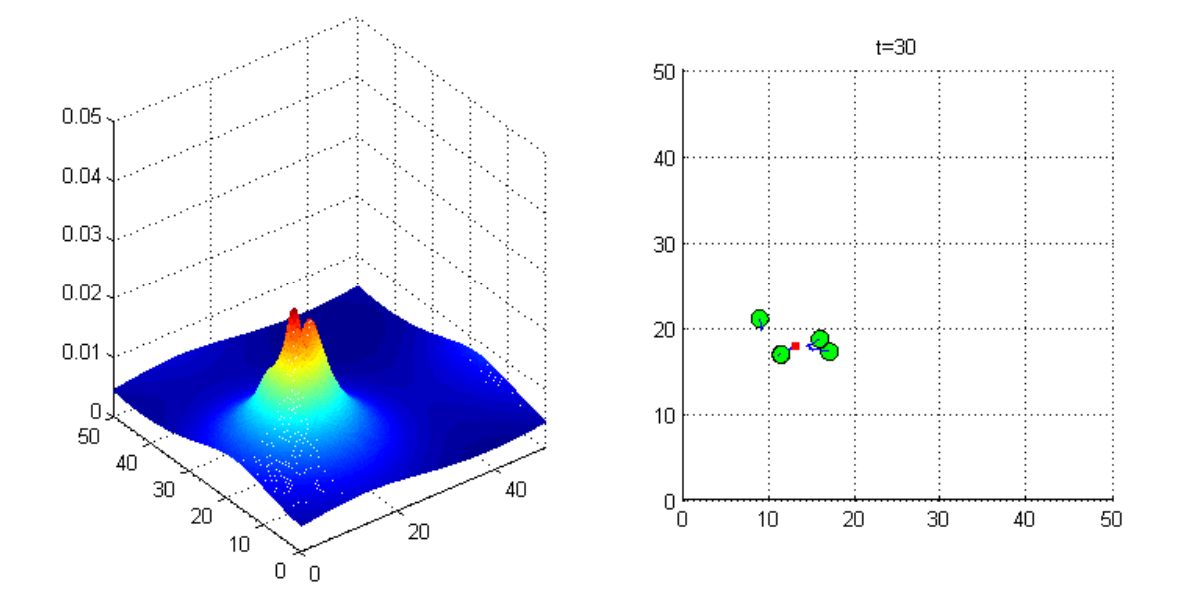}}
\subfigure[]{\includegraphics[width=0.45\textwidth]{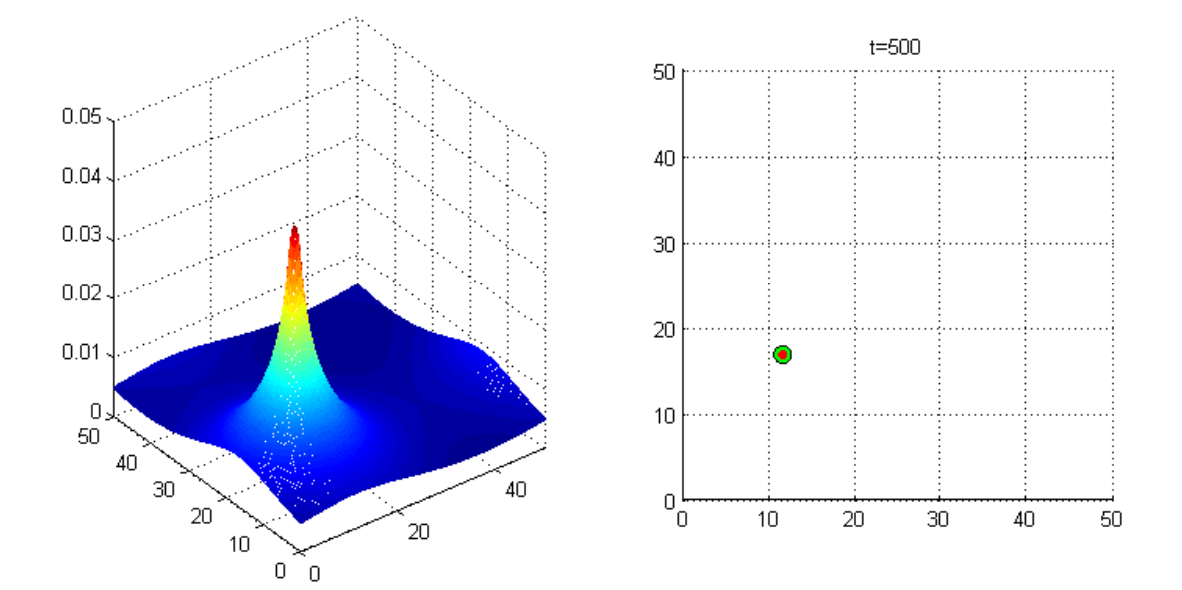}}
\caption{\textbf{Test 1}. Numerical simulation with parameters $\sigma=0.5$, $\beta=5$, $\gamma=2\times 10^2$, $D=2\times 10^2$, $\xi=0.5$, $V_{0,\max}=3$, and $\mathbf{X}_{0}$ randomly taken in the red square shown in (a) (Section \ref{sub-sec-tests}).}
\label{fig:test1}
\end{figure}

In the following numerical test we consider a spatial domain $\Omega=[0,50]\times[0,50]$ with periodic boundary conditions, and we choose a suitable time interval of observation $[0,T]$. For the initial data we fix $f(\mathbf{x},0)=0$ and, for $i=1,\dots,N$, $\mathbf{X}_i(0)=\mathbf{X}_{i0}$, $\mathbf{V}_i(0)=\mathbf{V}_{i0}$. In particular $\mathbf{X}_{i0}$ is chosen as a random vector, such that all the particles at $t=0$ are contained in a suitable initial region, fixed in the domain. Then $\mathbf{V}_{i0}=\left(V_{i0}\cos\theta_{i},V_{i0}\sin\theta_{i}\right)$ are chosen with $V_{i0}$ random numbers in $\left[0,V_{0,\max}\right]$, and $\theta_{i}$ random numbers in $[0,2\pi]$. 
We set the parameters $\sigma=0.5$, $\beta=5$, $\gamma=2\times10^2$, $D=2\times 10^2$, $\xi=0.5$, $V_{0,\max}=3$, and we consider $N=10$ particles located in $\mathbf{X}_{0}$ as in Figure \ref{fig:test1} (a). The time interval of observation is $[0,500]$. Spatial and temporal discretizations are given respectively by $\Delta x=\Delta y=0.125$ and $\Delta t=10^{-5}$.
\par
\begin{figure}[!ht]
\centering
\subfigure[]{\includegraphics[width=0.3\textwidth]{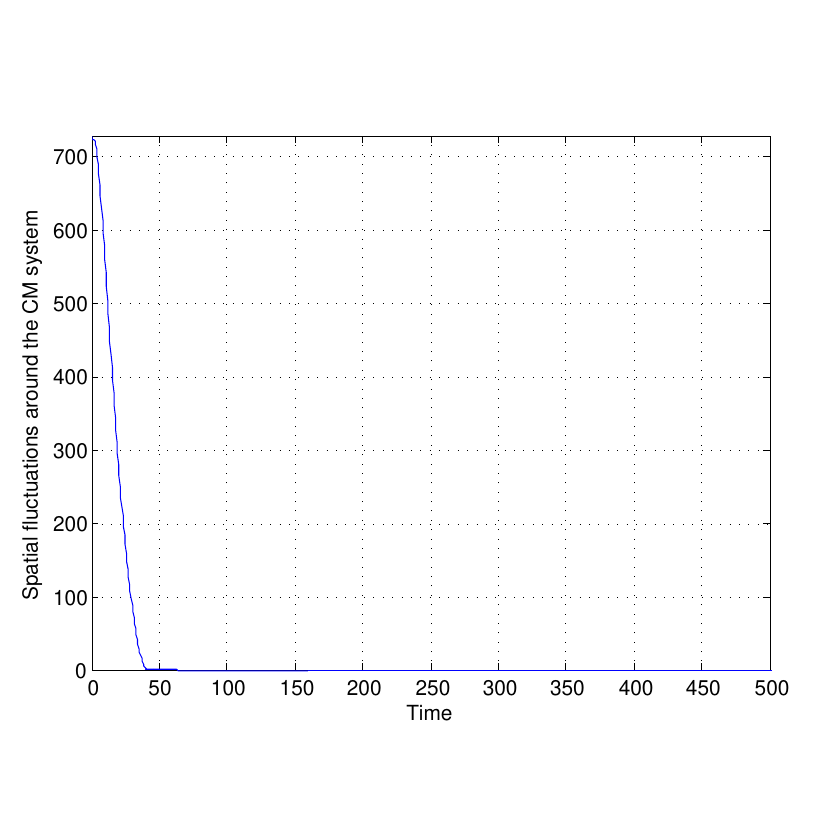}}
\hspace{-0.6 cm}
\subfigure[]{\includegraphics[width=0.3\textwidth]{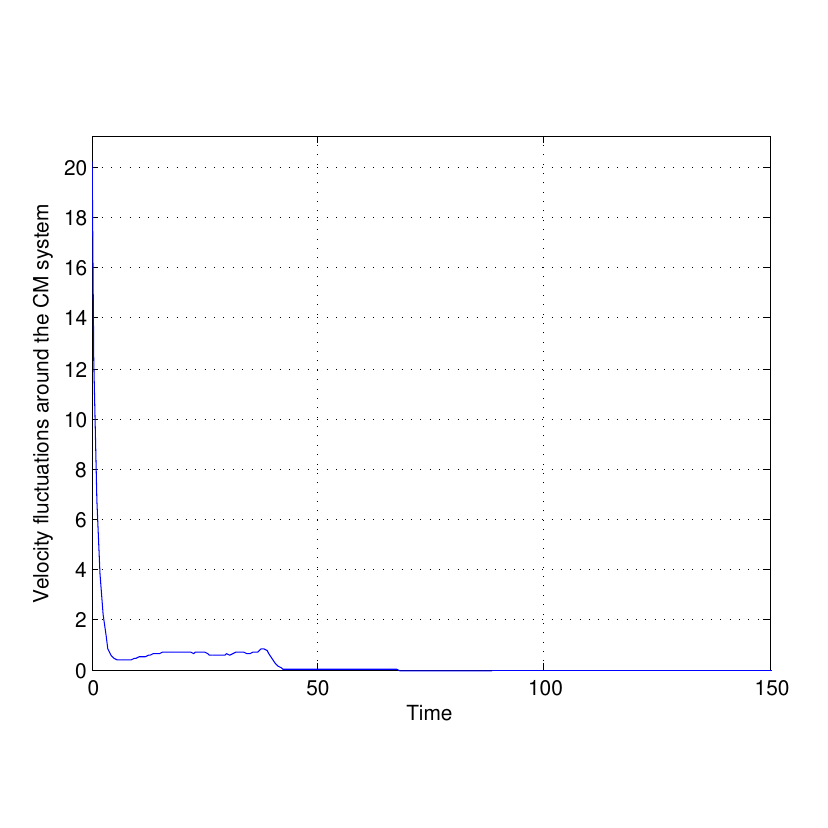}}
\vspace{0cm}
\hspace{-0.6 cm}
\subfigure[]{\includegraphics[width=0.3\textwidth]{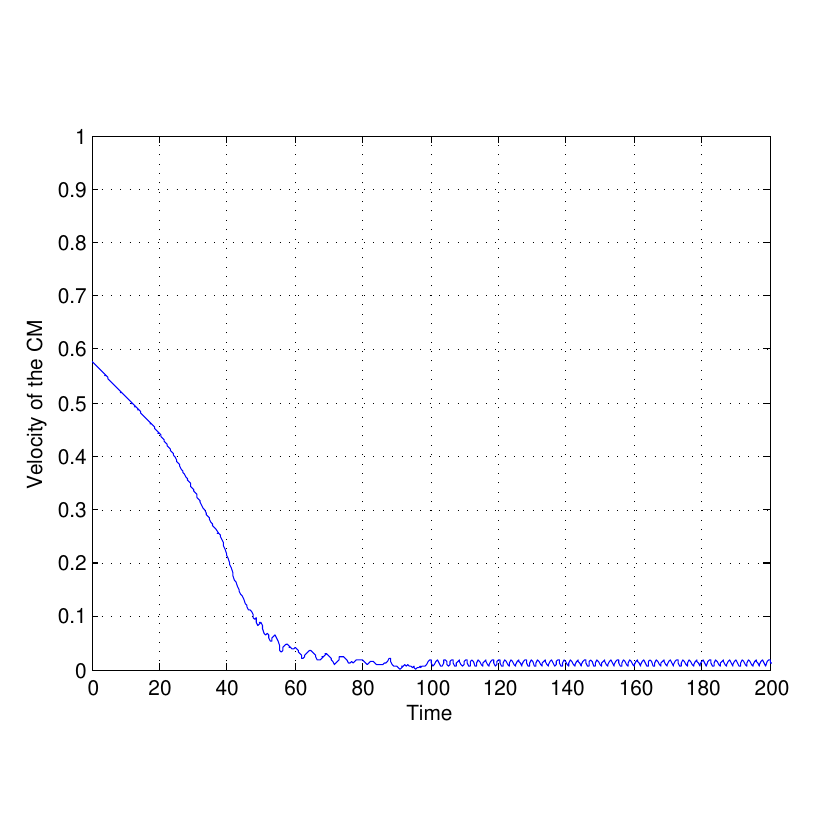}}
\caption{\textbf{Test 1}. (a)--(b) Functions $Fl_{X}(t)$ and $Fl_{V}(t)$ versus time (x-axis shows only a part of the time domain) as defined in Section \ref{sub-sec-tests}. (c) $\left\|\mathbf{V}_{\text{CM}}(t)\right\|$ versus time.}
\label{fig:test1-plot}
\end{figure}

 Figure \ref{fig:test1} shows four time steps of the numerical simulation: $t=0,5,30,500$. We plot on the left the chemoattractant concentration $f(\mathbf{x},t)$, while on the right the positions and the velocities of the particles in the spatial domain. The red square at $t=0$ is the region in which the initial positions are taken. The red marker indicates the centre of mass of the system, and the blue arrows are the velocity vectors. We observe an initial stage in which the particles tend to move somewhat aligned about until $t=5$ (Figure \ref{fig:test1} (b)), then they begin to converge to their centre of mass about at $t=30$ (Figure \ref{fig:test1} (c)), finally all particles stop in a same position (Figure \ref{fig:test1} (d)). In Figure \ref{fig:test1-plot} (a)--(b) are shown the spatial and velocity fluctuations around the centre of mass system
$$Fl_{X}(t):=\sum_{i=1}^N \left\|\mathbf{X}_{i}(t)-\mathbf{X}_{\text{CM}}(t)\right\|^2,\quad Fl_{V}(t):=\sum_{i=1}^N \left\|\mathbf{V}_{i}(t)-\mathbf{V}_{\text{CM}}(t)\right\|^2.$$
as a function of the time. For $t> 60$ $Fl_{X}(t)$ and $Fl_{V}(t)$ are less than $10^{-10}$. Notice that the square root of $Fl_{X}(t)$ and $Fl_{V}(t)$ is proportional to the standard deviations of $\mathbf{X}_i(t)$ and $\mathbf{V}_i(t)$ with respect to the position and velocity of the centre of mass. Figure \ref{fig:test1-plot} (c) displays the norm of the velocity of the centre of mass $\left\|\mathbf{V}_{\text{CM}}(t)\right\|$ versus time. For $t> 88$ this velocity is less than $2\times 10^{-2}$.
%
%
%
%
%
%
%
%
%


%
%
%
%
%
%
\section{Conclusions}\label{sec:conclusion-flock}
In this paper we have proposed an extension of the Cucker-Smale model introducing a mathematical model for collective motion driven by an alignment and a chemotaxis effect. We have adopted a hybrid description, discrete for the particles and for the motion equations, and continuous for the molecular level, containing the equation for the chemical signal.
\par We have studied our model by both an analytical and a numerical point of view. By the analytical point of view, using a fixed--point theorem, we have proved local and global existence and uniqueness of the solution of the nonlinear system. Moreover, we have investigated the asymptotic behaviour of the linearised system. We have proved the asymptotic convergence of the particles in their centre of mass with same velocity. Then the velocity of the centre of mass is proved to go time-asymptotically to zero. By a numerical point of view this property has been tested on the full nonlinear system, finding a complete concordance with the analytical results.
\par Future perspectives can concern the extension of the stability result, proved on the linearised system, to the full nonlinear model. Furthermore, other interactions could be taken into account in the collective motion, such as adhesion-repulsion, damping or lateral inhibition terms, similar to those introduced in \citet{dicostanzo}. This would be interesting in view of studying, in an analytical framework, the morphogenetic process arising in the lateral line development, and establish results in relation to the neuromast formation and deposition.
\bibliography{biblio-flock_chemo} 
\bibliographystyle{spbasic}
%
%
%
%
%
\end{document}